\documentclass[10pt, a4paper]{article}

\usepackage{amsmath,amssymb, amsthm}
\usepackage[utf8]{inputenc}
\usepackage{mathrsfs}
\usepackage{mathtools,setspace}
\usepackage{fullpage}
\usepackage[affil-it]{authblk}

\usepackage{verbatim}

\usepackage{hyperref}

\usepackage{xcolor}
\hypersetup{
  colorlinks   = true, 
  urlcolor     = {blue!90!black}, 
  linkcolor    = {blue!90!black}, 
  citecolor   = {red!90!black} 
}

\newcommand{\R}{\mathbb{R}}
\newcommand{\N}{\mathbb{N}}
\newcommand{\Z}{\mathbb{Z}}

\newcommand{\X}{\mathcal{X}}
\newcommand{\0}{\mathbf{0}}
\renewcommand{\d}{\mathrm{d}}

\theoremstyle{plain}
\newtheorem{thm}{Theorem}[section]
\newtheorem{prop}[thm]{Proposition}
\newtheorem{lem}[thm]{Lemma}

\theoremstyle{definition}

\newtheorem{defi}[thm]{Definition}
\newtheorem{rmk}[thm]{Remark}

\author{Viktor Bezborodov\thanks{Email: \texttt{vbezborodov@math.uni-bielefeld.de}} \quad
Yuri Kondratiev \thanks{Email: \texttt{yukondrat@yandex.ru}} \quad
Oleksandr Kutoviy \thanks{Email: \texttt{kutoviy@math.uni-bielefeld.de}} }

\affil{\emph{Bielefeld University, Faculty of Mathematics}}

\title{Lattice birth-and-death processes}

\begin{document}

\maketitle

\begin{abstract}
Lattice birth-and-death Markov dynamics
of particle systems with spins from $\Z _+$
are constructed as unique solutions to certain 
stochastic equations.
Pathwise uniqueness, strong existence, Markov property
and joint uniqueness in law are proven, and 
a martingale characterization of the process
is given. 
Sufficient conditions for the existence of 
an invariant distribution are formulated 
in terms of Lyapunov functions. 
We apply obtained results to discrete analogs
of the Bolker--Pacala--Dieckmann--Law model 
and an aggregation model.

\end{abstract}

\section{Introduction}

The evolution of a birth-and-death process admits the following
description. Two functions characterize the 
 development in time,
 the birth rate
 $b$
 and the death rate
 $d$. 
 If the system is in state $\eta \in \Z _+ ^{\Z ^\d}$
at time $t$,
then the probability that the number of points at a site $x \in \Z ^\d$ 
is increased by $1$ (``birth'') over the next time interval 
of length $ \Delta t$ is 
\[
 b(x,\eta)\Delta t + o(\Delta t),
\]
 the probability that the number of points 
at the site  $x$ is decreased by $1$ (``death'') over the next
time interval of length $\Delta t$ is 
\[
 d(x,\eta)\Delta t + o(\Delta t),
\]
and no two changes occur at the same time. Put differently,
a birth at the site $x$ occurs
at the rate
 $b(x,\eta)$, a death at the site $x$
 occurs at the rate $d(x,\eta)$, and no two events,
 births or deaths,
 happen simultaneously.
 
 The (informal) generator of such process is 
 
 \begin{equation} \label{the generator}
  L F (\eta) = \sum\limits _{x \in \Z ^\d} b(x, \eta) 
[F(\eta ^{+ x}) -F(\eta)] +
\sum\limits _{x \in \Z ^\d} d(x, \eta) [F(\eta ^{- x}) - F(\eta)],
 \end{equation}
where  
\[ \eta ^{+ x}(y) =
 \left\{
  \begin{array}{l l}
    \eta (y), & \text{if }  y \ne x, \\
    \eta (y)+1, & \text{if }  y = x,
  \end{array} \right. \ \ \ \
\eta ^{- x}(y) =
 \left\{
  \begin{array}{l l}
    \eta (y), & \text{if }  y \ne x,  \text{ or if } y=x, \eta (x)=0 \\
    \eta (y)-1, & \text{if }  y = x,  \eta (x) \ne 0.
  \end{array} \right.
\]

Birth-and-death processes we consider here 
correspond to lattice interacting particle systems
with a non-compact (single) spin space
and general transition rates.
The 
existence of the underlying stochastic dynamics
is not obvious. The first result of this
article is the construction 
of the corresponding Markov process.
Following ideas of Garcia and Kurtz \cite{GarciaKurtz},
we construct the process as a unique solution
to a certain stochastic integral equation
with Poisson noise. 

Birth-and-death processes constructed in \cite{GarciaKurtz}
are represented by a collection of points
in a separable complete metric space.
The scheme proposed there covers the case of $\Z _+ ^{\Z ^\d}$-valued processes.
The existence and uniqueness theorem, \cite[Theorem 2.13]{GarciaKurtz},
was obtained under the assumption that the death rate is constant,
which in the settings of this paper corresponds to 
$d(x, \eta ) = \eta (x)$ (although the existence
was shown under more general conditions). Kurtz and Protter 
\cite{KP96} give the uniqueness result
for stochastic equations driven by
semimartingale random measure (which include Poisson measures)
with in some sense Lipschitz coefficients.
 Theorem \ref{core thm} in this paper covers 
more general, not necessarily Lipschitz birth and death rates. 

A growing interest to the study of spatial birth-and-death processes 
which we have recently observed is stimulated by, among other things, an important role 
which these processes play in several applications. For example, in spatial plant ecology, 
a general approach to the so-called individual based models
was developed in a series of works, 
see e.g. \cite{BolkPacala, BolkPacala2, DieckmannLaw, MDL04, Ecology}
 and references therein. 
These models are represented by birth-and-death Markov processes in continuous configuration space
(over $\mathbb{R}^{d}$) with specific rates $b$ and $d$ which reflect biological notions
such as competition, establishment, fecundity etc. 
Other examples of birth-and-death processes may be found in mathematical physics,
see e.g. \cite{KKZh06, KKM10, Semigroupapproach} and references therein. 
Usually, lattice models can be compared with continuous ones
if we discretize the space $\mathbb{R}^{d}$ by partitioning it into cubes with centers at vertices of the lattice.
It is worth pointing out that in many applications
a lattice version can be constructed  which
will stochastically dominate the original continuous model.
Of course, such comparison arguments seem to be loose, because
the construction of the continuous original process is in general a very difficult problem.
 Nevertheless, 
 we may hope
  to deduce
 a priori information for the continuous process
from the corresponding lattice one.

There is an enormous amount of literature 
related to interacting particle systems in $\Z _+ ^{\Z ^\d}$.
Systems with a non-compact discrete spin space
appear as early as 1970 in the work of Spitzer \cite{Spit70},
where the invariant product measure
was constructed for the
zero range interaction.
The zero range process was constructed 
in a companion paper
by Holley \cite{Hol70} and later
by Andjel \cite{And82} under more general conditions
(see also Balázs et al. \cite{BRASS07}).
The zero range process represents 
the dynamics of hopping particles 
with the condition
that the jump rate depends only on
the number of particles 
at the departure site.
Various generalizations of 
the zero range process
have been considered, 
for example 
the so-called 
misanthrope process
was introduced by Cocozza-Thivent \cite{misanthrope85}. 
The zero range process 
has been extensively studied
ever since and 
has quite a few applications
in mathematical physics,
see e.g. the review by Evans and Hanney \cite{EH05}.
Another class of hopping particle models was 
considered for example by
Kesten and Sidoravicius \cite{KesS05}
(see also \cite{KesS08}),
where
an interacting particle system with 
non-trivial interaction was constructed and studied.
The system
models a spread of a rumor or infection 
and involves infinitely many particles.
We note that the birth-and-death systems
are of course different from the ones listed in this paragraph,
since the basic operations are the `addition' and `deletion'
of particles instead of the `replacement'.
Nor are
 the birth-and-death systems 
included in the class of
linear systems
(see e.g. \cite[Chapter 9]{Lig85} and references therein).

The scheme proposed by 
Etheridge and Kurtz \cite{EK14} covers
a wide range of systems
and applies to discrete and continuous models.
Their approach is based on, among other things, 
assigning a certain mark (`level')
to each particle and letting this mark evolve
according to a certain law. 
A critical event, such as a birth or death,
occurs when the level hits some threshold.
This scheme allows to consider 
multiple events and 
independent thinning, however
it seems to us that dynamics
with only a very specific types of interaction between particles
can be treated.
Penrose \cite{Pen08} gives a general existence result 
for particle systems with local interaction 
and uniformly bounded jump rates 
but non-compact spin space.
The results of \cite{Pen08} cannot be applied
to the systems discussed in the present paper 
since the rates are not supposed to be bounded.
Such systems are especially complicated for analysis.
For this reason 
the existence of the microscopic stochastic dynamics 
is sometimes simply assumed, see for example Bal{\'a}zs et al. \cite{Assumeexistence}.

The state space of our process will be
\[
 \X : = \Big\{ \eta \in \Z _+ ^{\Z ^\d} : \sum\limits _{x \in \Z ^\d}
 w (x) \eta (x) < \infty \Big\},
\]
where
$w$ is a summable positive even function.
Such subspaces of the product space 
naturally arise 
in the analysis of systems with unbounded transition 
rates because the process
started from arbitrary $\eta \in \Z _+ ^{\Z ^\d}$
need not exist; compare with
Liggett and Spitzer \cite{LS81}
and
\cite[(1.2)]{And82}.

In the present paper we have developed 
a technique which allows to construct 
lattice birth-and-death process
with unbounded transition rates.
We also mention that although our lattice
is given by $\Z ^\d$, the approach 
to construction that we use should
 work for an arbitrary connected bounded degree graph, 
 provided, of course, that the 
 assumptions are appropriately modified.
 Indeed,
 in our assumptions and proofs
 we use only the graph distance 
 on $\Z ^\d$ given by 
 \[
 |x-y|_1 := \sum\limits _{j=1 } ^\d |x_j - y_j|
 \]
 for $x = (x_1, ...,x_\d)$, $y = (y_1, ...,y_\d)$.

The paper is organized as follows. 
In Section \ref{results}
we collect the main results.
The first result of the paper is Theorem \ref{core thm}
which is an extension of 
the research done in the thesis \cite[Chapter 5]{Bezborodovthesis}.
A martingale characterization
of the constructed process and sufficient conditions
for the existence of an invariant distribution 
are given.
The proofs of the theorems from 
Section \ref{results} 
as well as some further comments
are given in
Sections \ref{the construction} through \ref{inv measure}.
In Section \ref{extinct surv} we  discuss survival 
of the process for a model with local death rate
and independent branching birth rate. We use comparison 
with the contact process to establish existence
of a critical value of the birth rate parameter.

\section{The set-up and main results}\label{results}

 Let $T>0$, and let 
 $N_1, N_2$ be 
 Poisson point processes on $ \Z ^\d \times \R_+ \times \R_+ $
 with 
intensity measure $\# \times ds \times du $, 
where
$\# $ is the counting measure on $ \Z ^\d$.
 Consider the equation
 
 \begin{equation}\label{se lat}
\begin{split}
\eta _t (x) = \int\limits _{(0,t] \times [0, \infty ) }
I _{ [0,b(x,\eta _{s-} )] } (u) N_1 (x,ds,du) & \\
- \int\limits _{(0,t] \times [0, \infty ) }
I _{ [0,d(x,\eta _{r-} )] } (v) N_2 (x,dr,dv) + \eta _0 (x), 
\end{split}
\end{equation}
where 
$(\eta _t)_{t\in[0,T]}$ is a c\`adl\`ag
 $\X$-valued solution process, 
 $x \in \Z ^\d $,  $ \eta _0$
is a (random) initial condition,
 $b,d$ are birth and death rates.
 We require
processes $N_1, N_2, \eta _0$ to be independent
of each other.
The first (second) term on the right hand side 
represents the number of births (deaths, respectively)
for $(\eta _t)$ at $x$ before $t$.
The integrals on the right-hand side are
taken in the Lebesgue—Stieltjes sense:
if for example 
$N _1 (\{ x\} \times \R _+ \times \R _+) = \sum _{i} \delta _{(x,s_i,u_i)}$,
then 
\[
 \int\limits _{(0,t] \times [0, \infty ) }
I _{ [0,b(x,\eta _{s-} )] } (u) N _1 (x,ds,du) = 
\sum\limits _{i: \ 0< s_i \leq t} 
I _{ [0,b(x,\eta _{s_{i}-} )] } (u_i).
\]
The theory of integration with respect to 
Poisson point processes can be found in 
Chapter 2 of \cite{IkedaWat}.
 Equation \eqref{se lat} is understood
in the sense that the equality
holds a.s. for all
$x \in \Z ^\d$ and $t \in (0,T]$.

We note here that equation \eqref{se lat}
is designed in such a way that 
the solution process 
satisfies the heuristic description 
at the beginning of the introduction.
In Proposition \ref{poignant}
we will see a formal connection of a unique solution
to \eqref{se lat} and the heuristic generator 
given in \eqref{the generator}.

\emph{Assumptions on} $\eta _0$, $w$, $b$ \emph{and} $d$.
Let us fix the assumptions we use throughout the paper.
Let $w$ be a summable positive even function, $\sum\limits _{x \in \Z ^\d} w (x) < \infty$.
Denote by $\X$ the set
\[
 \{ \eta \in \Z _+ ^{\Z ^\d} : \sum\limits _{x \in \Z ^\d}
 w (x) \eta (x) < \infty \}.
\]

We equip $\X$ with the topology induced by the distance
\[
 d_{\X} (\eta ,\zeta ) = \sum\limits _{x \in \Z ^\d} w (x) |\eta (x) - \zeta (x)|.
\]
Note that $(\X, d_{\X})$ is a complete separable metric space and that
convergence in
 $\X$ implies pointwise convergence: 
 if $\eta _k \in \X$, $k=0,1,...$, $\eta \in \X$  and $\eta _k \to \eta $ in $\X$,
 then $\eta _k (x) \to \eta (x)$ for any $x \in \Z ^\d$ and,
 since $\eta _k(x)$ is a natural number or zero, 
 $\eta _k (x) = \eta (x)$ for all but finitely many $k$.
 We require that 
 \[
 E \sum\limits _{x \in \Z ^\d}
 w (x) \eta _0 (x) < \infty.
\]
Clearly, the latter implies that a.s. $\eta _0 \in \X$.

The birth and death rates $b$ and $d$ are functions 
defined on $\Z ^\d \times \X$ and taking values in $\R _+$.
 We assume throughout that the following conditions are satisfied:
 \begin{equation}\label{bigot}
  \begin{split} 
   \text{if } \quad \xi, \eta \in \X, \  x \in \Z ^\d \text{ and } \ \xi (x) \geq \eta (x), \quad \text{then } \\
   b(x, \xi) - b(x, \eta ) \leq \sum\limits _{y \in \Z ^\d} a(x-y) |\xi (y) - \eta (y)|,
  \end{split}
 \end{equation}
  \begin{equation}\label{spaz}
  \begin{split} 
   \text{if } \quad \xi, \eta \in \X, \ x \in \Z ^\d \text{ and } \ \xi (x) \geq \eta (x), \quad \text{then } \\
   d(x, \xi) - d(x, \eta ) \geq - \sum\limits _{y \in \Z ^\d} a(x-y) |\xi (y) - \eta (y)|,
  \end{split}
 \end{equation}
 and
 \begin{equation}\label{foible}
  \sum\limits _{y \in \Z ^\d} w (y) a (x-y) \leq C_{w,a} w (x),  \quad \quad
  x \in \Z ^\d
 \end{equation}
where $a:  \Z ^\d \to \R _+$ is a summable even function,
$C_{w,a} >0$.
Denote by $\0$ the ``zero'' element of $\X$:
  $\0 \in  \Z _+ ^{\Z ^\d}$, $\0 (x) = 0$, $x \in \Z ^\d$.
  If there are no particles at a site then no death can occur,
  so $d$ should satisfy
 \begin{equation*}
   d(x,\eta) =0, \quad \text{ whenever } \  \eta (x) = 0.
 \end{equation*}
We also require
\begin{equation}
  \sum\limits _{x \in \Z ^\d} w (x) b(x, \0) < \infty.
  \end{equation}

  For some possible choices of $a$ and $w$ satisfying these conditions and for a few
  examples, see Remark \ref{examples} below.

  We say that a Poisson point process $N$ on $\Z ^\d \times \R _+ \times \R _+$
  is \emph{compatible} with a right-continuous complete filtration $\{ \mathscr {F} _t \}$
  if all random variables of the form $N(\{ x \} \times [a,b] \times U )$, 
  $x \in \Z ^\d$, $0\leq a < b \leq t$, $U \in \mathscr{B} (\R _+)$,
  are $ \mathscr {F} _t $-measurable, and, in addition, 
  all random variables of the form $N(\{ x \} \times [a,b] \times U )$,
  $x \in \Z ^\d$, $t < a < b $, $U \in \mathscr{B} (\R _+)$,
  are independent of $ \mathscr {F} _t $.

\begin{defi} \label{weak solution lat}
 A \emph{(weak) solution} of equation \eqref{se lat} is a triple 
 $(( \eta _t )_{t\in[0,T]} , N_1 , N_2)$, 
 $(\Omega , \mathscr{F}, P) $,
 $(\{ \mathscr {F} _t  \} _ {t\in[0,T]}) $, where 

  (i) $(\Omega , \mathscr{F} , P) $
  is a probability space, 
 $\{ \mathscr {F} _t  \} _ {t\in[0,T]}$ is an increasing, 
right-continuous
 and complete filtration of sub-$\sigma$-algebras of $\mathscr {F}$,

  (ii)  $( \eta _t )_{t\in[0,T]} $ is a c\`adl\`ag 
  process in 
$\X$,
adapted to $\{ \mathscr {F} _t  \} _ {t\in[0,T]} $, such that
$${E\sum\limits _{x \in \Z ^\d} w (x) \sup\limits _{t \in [0,T]}  \eta _t (x) < \infty},$$
  (iii)
$N_1 , N_2$ are independent
 Poisson point processes with measure intensity  
$\# \times ds \times du $, compatible with 
$\{ \mathscr {F} _t  \} _ {t\in[0,T]} $,
  
  (iv) all integrals in \eqref{se lat} are well-defined, and $E \int\limits _0 ^T
  [b(x,\eta _{s-}) + d(x,\eta _{s-})] ds < \infty$ for every $x \in \Z ^\d$.

  (v) equality \eqref{se lat} holds a.s. for all ${t\in[0,T]}$
  and $x \in \Z ^\d$.

\end{defi}

\begin{rmk}

The definition above as well as many of the definitions
and theorems
below 
can be extended to the case of the time interval $[0,\infty)$
in an obvious manner.
\end{rmk}

  \begin{defi}\label{strong solution lattice}
  A solution is called \emph{strong} if $( \eta _t )_{t\in[0,T]} $
  is adapted to the completion under $P$ of the filtration
  $$
  \mathscr{S} _t = \sigma \{ \eta_0 , N_k(\{x\} \times [0,q] \times C), 
x\in \Z ^\d, C\in \mathscr{B} (\R_+), q\in [0,t] , k=1,2 \}.
$$

 \end{defi}

 For complete $\sigma$-algebras $\mathscr{A}_1$ and $\mathscr{A}_2$,
 let $\mathscr{A}_1 \vee \mathscr{A}_2$ 
 be the smallest complete $\sigma$-algebra 
 containing both $\mathscr{A}_1$ and $\mathscr{A}_2$.
 
  \begin{defi}\label{path uniq def}
 We say that \emph{pathwise uniqueness} holds for equation \eqref{se lat} and an
initial distribution $ \nu $ if, whenever the triples 
$(( \eta _t )_{t\in [0,T]} , N_1 ,  N _2)$, $(\Omega , \mathscr{F} , P) $,
 $(\{ \mathscr {F} _t  \} _ {t\in [0,T]}) $ and 
 $(( \bar \eta _t )_{t\in [0,T]} , N_1 ,  N _2)$, $(\Omega , \bar {\mathscr{F}} , P) $,
 $(\{ \bar{\mathscr {F}} _t  \} _ {t\in [0,T]}) $ are weak solutions of \eqref{se lat} with 
$P \{ \eta _0 = \bar{\eta} _0 \} = 1 $, $Law (\eta _0) = \nu $, 
and such that the processes $N_1$, $N_2$ are compatible 
with  $ \left\{ \mathscr {F} _t \vee \bar{\mathscr {F}} _t  \right\} _ {t\in [0,T]}  $,
we have $P \{ \eta _t = \bar{\eta} _t , t \in [0,T] \} = 1$
(that is, the processes $(\eta _t) , (\bar{\eta _t}) $ are indistinguishable).

\end{defi}

\begin{defi}
 We say that \emph{joint uniqueness in law} holds for \eqref{se lat} 
 if the triple $(( \eta _t )_{t\in[0,T]} , N_1 , N_2)$
 has the same distribution in 
 $D _{\X}[0,T] \times D _{\mathbf{E}} [0,T] \times D _{\mathbf{E}} [0,T]$ for every weak solution (cf., e.g., 
 \cite[Definition 2.9]{KurtzYamWatan}). 
\end{defi}
 Here $D _{\X}$ and $D _{\mathbf{E}}$ are 
 the spaces of c\`adl\`ag paths 
 over the corresponding spaces
 equipped with the Skorokhod topology,
 and $\mathbf{E}$ is the space of locally finite
 simple counting measures 
 on $\Z ^\d \times \R _+$ with the minimal $\sigma$-algebra such that every set of the form
 \[
  \{ \gamma \in \mathbf{E} \mid \gamma(Q) \in B \}, 
  \quad \quad Q \in \mathscr{B}(\Z ^\d \times \R _+), \ B \in \mathscr{B} (\R _+)
 \]
is measurable, and endowed with the metric compatible
with the vague topology (also called the \emph{space of locally finite
configurations}; see e.g. \cite[Appendix A2]{Kallenbergfound} or \cite{KondKuna},
and references therein).

Now we formulate the existence and uniqueness theorem, which will be proven in the next section.

\begin{thm} \label{core thm}
Under the above assumptions 
 pathwise uniqueness, strong existence and uniqueness in law hold for 
 equation \eqref{se lat}. The unique solution 
 is a Markov process.
\end{thm}

   \begin{rmk} \label{examples}
   The conditions on $b$, $d$ and $\eta _0$ given
   in terms of functions $w$ and $a$
   may seem somewhat indirect.
   In fact, given $a$ satisfying our assumptions,
   it is always possible to construct $w$
   satisfying our assumptions following the scheme
   in \cite{LS81} (to ensure that $w$ is even,
   we should take the function $\beta$ there to be even).
   Here we point out three
   possible choices of $w$ and $a$. 
   
   (i) $w(x) = e ^{- q |x|_1}$, $a(x) = c e ^{- p |x|_1}$,
   $p > q >0$, $c>0$;
   
   (ii) $w(x) = e ^{- q |x|_1}$, $a(x) = c I \{ |x|_1 \leq k\}$,
   $q >0$, $c>0$, $k \in \N$ ;
   
   (iii) $w(x) = \frac{1}{1 + |x|_1 ^{-q}}$, $a(x) = \frac{c}{1 + |x|_1 ^{-p}}$,
   $p > q > \mathrm{d}$, $c>0$.
  
  \end{rmk}
 
 Now we give two examples 
 where Theorem \ref{core thm} applies, and we can obtain
 discrete counterparts of some
 continuous particle systems as 
 unique solutions to \eqref{se lat}.
 
 \emph{A discrete version of the Bolker--Pacala--Dieckmann--Law model}
 also known as spatial stochastic logistic model 
  (Bolker and Pacala \cite{BolkPacala, BolkPacala2}; Dieckmann and Law et al. 
  \cite{DieckmannLaw, MDL04}).
An individual-based description of this model is as follows:

(1) Existing individuals produce offsprings at a per capita fecundity rate. 

(2) A newly produced offspring is distributed (instantaneously) according to a dispersal kernel,
and it is assumed to establish (instantaneously) as a newborn
individual, which matures (instantaneously) and starts to produce offsprings. 

(3) Existing individuals may die for two reasons. Firstly, there is a constant background per capita mortality rate $m$, 
yielding an exponentially distributed
lifetime. Secondly, mortality has a density dependent component (self-thinning), 
so that competition among the individuals may also lead to death. 
The density dependent
component of the death rate of a focal individual is a sum of contributions
from all the other individuals within the entire $\mathbb{Z}^d$,
but the strength of the competitive
effect decreases with distance. 

The model is defined by
 \begin{equation} \label{flabbergast}
  b(x, \eta) = b_0 + \sum\limits _{y \in \Z ^\d} a_+ (x-y) \eta (y), 
 \end{equation}
 and 
  \begin{equation}\label{foolhardy}
  d(x, \eta) =  m_1 \eta (x) +  m_2 \eta (x)(\eta (x) - 1) + \sum\limits _{y \in \Z ^\d} a_- (x-y) \eta (y).
 \end{equation}
 In this model 
$(\eta _t)$ represents an evolution of a biological 
population with independent branching given 
by the kernel $a _+$, immigration at a constant rate $b_0$,
constant ``intrinsic''
mortality rate $m_1$, local competition rate $m_2$, and the competition kernel $a_-$.
The functions $a_+$ and $a_-$ are assumed to be summable.
The model with rates \eqref{flabbergast}
and \eqref{foolhardy} with $b_0 = 0$
can be regarded as a translation invariant discretization 
of the  Bolker--Pacala--Dieckmann--Law model 
studied 
by Fournier and M{\'e}l{\'e}ard
\cite{FournierMeleard}.

The stepping stone and 
superprocess versions
of the Bolker--Pacala--Dieckmann--Law model were 
considered by Etheridge \cite{Eth04}.
The process with a finite number of
particles in continuum was studied in
\cite{FournierMeleard}, where, among other things,
it was shown that
the superprocess version can be obtained as a scaling limit
of continuous processes.
Statistical dynamics were considered by 
Finkilstein et al. \cite{FKK09, FKKBPDL};
see also Ovaskainen et al. \cite{Ecology}.
 Unlike in the continuous model, in the discrete model
 we allow the particles to be at the same 
 place; otherwise the density would be bounded.
 
 \emph{A discrete version of an aggregation model}.
 Here the birth rate is either as in \eqref{flabbergast},
 or is given by a constant, $b(x,\eta) \equiv c >0$.
 The death rate is given by 
 \begin{equation*}
  d(x,\eta) = e ^{-c \sum\limits _{y \in \Z ^\d} \varphi (x-y) \eta (y)},
 \end{equation*}
or 
  \begin{equation*}
  d(x,\eta) = \frac{1}{1 + c \sum\limits _{y \in \Z ^\d} \varphi (x-y) \eta (y)} ,
 \end{equation*}
 where $c>0$ and $\varphi : \Z ^\d \to \R _+$.
 For a statistical dynamics corresponding to this model in continuum, 
 see \cite{Aggreg} and references therein;
 see also \cite{Bezborodovthesis} for continuous systems 
 with a finite number of particles.

 The following two propositions establish a rigorous relation 
 between the unique solution to \eqref{se lat} 
 and $L$ defined by \eqref{the generator}.
 To formulate the first of them, 
 let us consider the class $\mathscr{C} _b$ of cylindrical functions $F :
\mathcal{X} \to \R _+$ with bounded increments.
We say that $F$ has bounded increments if 
\[
\sup\limits _{\eta \in \mathcal{X}, x \in \Z ^\d}
\big| F(\eta ^{+x} ) - F(\eta ) \big| < \infty.
\]
We say that $F$ is cylindrical 
if $F ( \eta)$ depends on values of $\eta $ 
in finitely many sites only, i.e. for some $R=R_F>0$
 \begin{equation*}
   F (\eta) = F(\zeta) \text{ whenever } \eta (x) = \zeta (x) \text{ for all }
  x, |x|_1 \leq R.
 \end{equation*}
 We recall that
 $|x|_1 = \sum\limits _{j=1} ^\d |x_j|$ for $x= (x_1,...,x_d)$
  and that the filtration $\{ \mathscr{S} _t, t \geq 0\}$
 appeared in Definition \ref{strong solution lattice}.

\begin{prop} \label{poignant}
 Let $(\eta _t) _{t \geq 0}$ be a weak solution to \eqref{se lat}. Then 
 for any $F \in \mathscr{C} _b$
 the process
 \begin{equation} \label{hooky}
  F(\eta _t) - \int\limits _0 ^t L F (\eta _{s-})ds 
 \end{equation}
is an $\{ \mathscr{S} _t \}$-martingale.
In particular, the integral in \eqref{hooky} is a.s. well-defined.
\end{prop}

The next proposition says that
under some additional assumptions
the converse is true.

 \begin{prop} \label{filibuster}
  
  Let $(\eta _t ) _{t \in [0,T]}$ be a
  c\`adl\`ag $\mathcal{X}$-valued process
  defined on a probability space $(\Omega , \mathscr{F}, P) $,
  adapted to a right-continuous complete filtration
 $(\{ \mathscr {F} _t  \} _ {t\in[0,T]}) $
 and
  satisfying
  \begin{equation*}
  \eta _t = \eta _0 + \sum\limits _{0 < s \leq t} (\eta _{s} - \eta_{s-}),  t \in [0,T],
  \end{equation*}
 and
  \begin{equation*}
     \sup\limits _{t \geq 0} \sum\limits _{x \in \Z ^\d} 
     |\eta _{t} (x) - \eta _{t-} (x)| \leq 1 
  \end{equation*}
a.s.
Assume that the probability 
space $(\Omega , \mathscr{F}, P) $ and 
  the filtration
 $(\{ \mathscr {F} _t  \} _ {t\in[0,T]}) $
 are rich enough to support 
 required randomization processes,
 and that 
\begin{equation}
 \eta _{t} (x) = \eta _{0} (x) + \eta ^{(b)} _{t} (x) -\eta ^{(d)} _{t} (x),
\end{equation}
where $ \{ \eta ^{(b)} _{t}, t\geq 0 \}$ and $ \{ \eta ^{(d)} _{t}, t\geq 0 \}$ 
are c\`adl\`ag non-decreasing $\mathcal{X}$-valued processes,
${\eta ^{(b)} _{0} = \eta ^{(d)} _{0} = \0 }$, 
such that for every finite $E_1, E _2 \subset \Z ^\d$
the process
\[
 \sum\limits _{x \in E _1} \eta ^{(b)} _{t} (x) + 
 \sum\limits _{y \in E _2} \eta ^{(d)} _{t} (y)
 - \int\limits _0 ^t  \sum\limits _{x \in E _1} b(x, \eta _{s-})ds
 - \int\limits _0 ^t  \sum\limits _{y \in E _2} d(y, \eta _{s-})ds
\]
is an $(\{ \mathscr {F} _t  \} _ {t\in[0,T]})$-martingale.
Then 
there exist independent Poisson point processes
$N_1$ and $N_2$ such that
the triple $(( \eta _t )_{t\in[0,T]} , N_1 , N_2)$, 
 $( \Omega , {\mathscr{F}},  P) $,
 $(\{ {\mathscr {F}} _t  \} _ {t\in[0,T]}) $
 is a weak solution to \eqref{se lat}. 
 \end{prop}

 \begin{rmk} \label{spin systems}
 Some classic
  interacting particle systems, including 
  the stochastic Ising model, the contact process
  and the voter model (see, e.g., \cite{Lig85})
  can be constructed using the above results.
  For example, the unique solution
  of \eqref{se lat}
  with initial condition
  $\beta \in \{0, 1 \}^{\Z ^\d} $,
  the death rate $d(x,\eta) = I \{ \eta (x) >0 \}$
  and the birth rate 
  \begin{equation} \label{smug}
 b (x,\eta ) =  b _{cont}(x,\eta ) := I_{ \{ \eta (x) =0  \} } \lambda \sum\limits
 _{y: |y-x|\leq 1} \eta (y)
 \end{equation}
  is the contact process with parameter $\lambda >0$
  and initial state $\beta$.
  This follows from the uniqueness of 
  solutions for the associated martingale problem, 
  see \cite[Theorem (4.12)]{HolStr76},
  and Proposition \ref{poignant}.
  
 \end{rmk}

 Sometimes we will denote the solution to \eqref{se lat}
 with initial condition $\eta _0 \equiv \alpha$, $\alpha \in \X$,
 by $(\eta (\alpha, t))_{t \in [0,T]}$, emphasizing the dependence on $\alpha$.

For $\alpha \in \mathcal{X}$ we denote by 
 $P _{\alpha}$ the law of $(\eta (\alpha, t))$,
 \[
  P _{\alpha} (H) 
  = P \{ (\eta (\alpha, t))_{t \geq 0} \in H \}
 \]
 for a measurable $H \in D _{\mathcal{X}} [0,\infty)$.

Note that $P _{\alpha}$ is well defined 
by Theorem \ref{core thm}.
Also, 
for every 
$H \in D _{\mathcal{X}} [0,\infty)$,
$P _{\alpha} (H) $
can be shown to 
be measurable in $\alpha$.
 
 Let $C_b (\X)$ be the 
 space of bounded continuous functions
 on $\X$ equipped
 with the supremum norm.
 For $\alpha \in \X$ and $f \in C _b(\X )$ we define
 
 \begin{equation}
  P ^t f (\alpha) : =  E f(\eta(\alpha, t))  \quad 
(= E _{\alpha} f(\eta _t)).
 \end{equation}

 The function $P ^t f$ is 
 continuous on $\X$. Indeed, by Lemma \ref{bloviate}
 and Gronwall's inequality
 \begin{equation}
 E  \sum\limits _{x \in \Z ^\d}  w (x) |\eta (\alpha, t) (x) - \eta (\beta, t) (x)|
 \leq
     E  \sum\limits _{x \in \Z ^\d}  w (x) |\alpha (x) - \beta (x)|
  \exp \{4 C_{w,a} t\},
 \end{equation}
 and hence by Lebesgue's dominated convergence theorem 
 $P ^t f$ is continuous.
 Therefore $P ^t$ is a bounded operator on $C_b (\X)$.

 A probability measure $\pi$ on $\mathcal{X}$ 
is called invariant for equation \eqref{se lat}
 if
\[
 \int P ^t f (\alpha) \pi (d \alpha)
 = \int f (\alpha) \pi (d \alpha)
\]
 for every $f \in C _b (\X)$.

Consider the following additional assumption:
there exists an even summable function \\
 ${v : \Z ^\d \to (0,\infty)}$ and a constant $C_{v,a} > 0$
 satisfying
 \begin{equation} \label{tipsy}
  \frac{v(x)}{w(x)} \to \infty, \quad x \to \infty,
 \end{equation}
 and
 \begin{equation} \label{salvo}
  \sum\limits _{y \in \Z ^\d} v (y) a (x-y) \leq C_{v,a} v(x), \quad x \in \Z ^d.
 \end{equation}

 Let $\X _v : = \{ \eta \in \X : \sum\limits _{x \in \Z ^\d} v(x) \eta (x) < \infty \}$
 and let $V:\X _v \to \R _+ $
 be given by 
 \[
 V (\eta) := \sum\limits _{x \in \Z ^\d} v(x) \eta (x), \quad \eta  \in \X _v.
 \]

\begin{thm} \label{schmooze}
  
Assume that there exists
an even summable function
$v: \Z ^\d \to (0,\infty) $ 
 such 
that \eqref{tipsy} and \eqref{salvo}
hold. 
 Also, assume that 
 for some constants $c_1, c_2 >0$
  \begin{equation}\label{clairvoyance}
  LV (\eta) \leq c_1 - c_2 V (\eta), \quad \textnormal{ for all } \eta \in \X _v.
 \end{equation}
Then there exists an invariant measure for equation \eqref{se lat}.
 
  \end{thm}
  
  Note that, for $\eta \in \X _v$, $LV (\eta)$
  may be equal to $- \infty$, in which case
   \eqref{clairvoyance} is fulfilled. 
 Using the theorem above, we can establish existence 
 of an invariant measure for 
 the model given by \eqref{flabbergast} and 
 \eqref{foolhardy}.

 \begin{prop} \label{beseech}
 Let $w(x) = e ^{-|x|_1}$.
Assume that $a_+ $ and $a_-$ in \eqref{flabbergast}
and \eqref{foolhardy} have the finite range property:
there exists $R>0$ such that $a_+(x) = a_-(x) = 0 $
whenever $|x|_1 \geq R$, and that $m_1, m_2 >0$. Then 
  equation \eqref{se lat} has
  an invariant measure.  
 \end{prop}

\section{Proof of Theorem \ref{core thm}}\label{the construction}

The statement of Theorem \ref{core thm}
is contained in 
Propositions \ref{un}, \ref{ex}, \ref{uniqueness in law thm} and \ref{Markov property},
which we prove below.

We start with the following Lemma.

\begin{lem} \label{brisk}
 For every $x\in \Z ^\d$
the maps
\begin{align*}
 \X \ni \xi \mapsto b(x, \xi)   \in \R _+, \notag
\\ 
 \X \ni \xi \mapsto d(x, \xi)   \in \R _+ \notag
\end{align*}
are continuous.
\end{lem}
\textbf{Proof}.
We give the proof for $ b$ only,
as the proof for $d$ 
can be done in the same way.
Fix $x\in \Z ^\d$ and $\xi \in \X$. 
Take $\delta \in (0, w (x))$
and $\eta \in \X$
such that 
$d _{\X} (\xi ,\eta) \leq \delta$. Then
$\xi (x) = \eta (x)$. We have by \eqref{bigot}
\[
  b(x, \xi) -  b(x, \eta ) \leq \sum\limits _{y \in \Z ^\d} a(x-y) |\xi (y) - \eta (y)|
\]
and
\[
  b(x, \eta) -  b(x, \xi ) \leq \sum\limits _{y \in \Z ^\d} a(x-y) |\xi (y) - \eta (y)|
\]
Hence
\[
 | b(x, \xi) - b(x, \eta )| \leq \sum\limits _{y \in \Z ^\d} a(x-y) |\xi (y) - \eta (y)|.
\]
Now, 
\eqref{foible} implies $w (y) a(x-y) \leq  C_{w,a} w (x)$, 
or, after swapping $x$ and $y$, ${w (x) a(x-y) \leq  C_{w,a} w (y)}$
and $a(x-y) \leq \frac{C_{w,a}}{w (x)} w (y)$.
Thus,
\[
 | b(x, \xi) -  b(x, \eta )| \leq 
 \frac{C_{w,a}}{w (x)} 
 \sum\limits _{y \in \Z ^\d} w (y) |\xi (y) - \eta (y)|. 
\]
\qed

Before treating equation \eqref{se lat}
in a general form,
let us consider the case of a ``finite'' initial condition. 
We call $\eta _0$ satisfying 
\begin{equation} \label{chafe}
 \sum _x \eta  _0 (x) < \infty \quad \text{ a.s. }
\end{equation}
and 
\begin{equation}\label{jettison}
 E \sum _x \eta  _0 (x) < \infty.
\end{equation}
a \emph{finite initial condition}.
Of course, \eqref{chafe} follows from \eqref{jettison}.

\begin{prop}\label{finite in cond}
Assume that there exist $c_1, c_2  \geq 0$ such that
\begin{equation} \label{smolder}
 \sum _{x \in \Z ^\d} b(x ,\eta ) \leq c_1 \sum _{x \in \Z ^\d} \eta (x) +c_2.
\end{equation}
 Then pathwise uniqueness and strong existence hold for
 \eqref{se lat} with a finite initial condition. 
 Furthermore, the unique solution $(\eta _t)_{t \in [0,T]}$ satisfies
\begin{equation}\label{scrounge-1}
   E \sum _x  \sup\limits _{t \in [0,T]} \eta _t (x) < \infty.
\end{equation}

\end{prop}

The proof can be done constructively, ``from one jump to another'',
following the proof of the existence and uniqueness theorem
for a similar equation
in continuous space settings, see \cite{Bez15} or \cite[Theorem 2.1.6]{Bezborodovthesis}.
The assertion \eqref{scrounge-1} follows from \eqref{smolder} and comparison with 
the Yule process.
The Yule process $(Z_t)_{t\geq0}$ is an $\N$-valued birth process
such that for all $n \in \N$
 \[
   P\{Z_{t+ \Delta t} - Z_t = 1 \mid Z_t = n \} = \mu n + o(\Delta t).
 \]
for some $\mu >0$; see e.g. \cite[Chapter 3]{Branch1} or \cite{Yuleproc},
and references therein.

Note that \eqref{scrounge-1} implies 
\begin{equation}\label{scrounge}
   E \sum _x w(x) \sup\limits _{t \in [0,T]} \eta _t (x) < \infty, \ \ \ t \in [0,T],
\end{equation}
since $w$ is summable and therefore bounded.

Consider now two solutions $(\eta ^{(k)} _t )$, $k =1,2$, to the equations
 \begin{equation}\label{flippant}
\begin{split}
\eta _t (x) = \int\limits _{(0,t] \times [0, \infty ) }
I _{ [0, b _k(x,\eta _{s-} )] } (u) N_1 (x,ds,du) \\
- \int\limits _{(0,t] \times [0, \infty ) }
I _{ [0, d _k (x,\eta _{r-} )] } (v) N_2 (x,dr,dv) + \eta ^{(k)} _0 (x)
\end{split}
\end{equation}
with finite initial conditions.

\begin{prop}\label{comparison finite}
Let $\eta ^{(1)}_0 $ and $ \eta ^{(2)}_0$ be finite initial conditions.
  Assume that almost surely $\eta ^{(1)}_0 \leq \eta ^{(2)}_0$, and
  
  (i) for any $\xi ^{(1)}, \xi ^{(2)} \in \X$ such that  $\xi ^{(1)} \leq \xi ^{(2)}$ and
  $\sum _{y \in \Z ^\d} \xi ^{(2)} (y) < \infty$,
\begin{equation} \label{culpable}
b_1(x, \xi ^{(1)}) \leq b_2 (x, \xi ^{(2)}), \ \ \ x \in \Z^\d,
\end{equation}

 (ii) for any $ x \in \Z^\d$ and $\xi ^{(1)}, \xi ^{(2)} \in \X$ such that $\xi ^{(1)} \leq \xi ^{(2)}$, 
  $\sum _{y \in \Z ^\d} \xi ^{(2)} (y) < \infty$ and  $\xi ^{(1)} (x) = \xi ^{(2)} (x)$,
\begin{equation}
 d_1(x,\xi ^{(1)}) \geq d_2 (x, \xi^{(2)}).
\end{equation}
Then 
\begin{equation} \label{culprit}
 \eta ^{(1)} _t \leq \eta ^{(2)} _t, \ \ \ t \in [0,T].
\end{equation}
Furthermore, the inclusion
\[
 \left\{ (t,x):  (\eta ^{(1)} _t (x) - \eta ^{(1)} _{t-} (x) )>0  \right\}
 \subset
 \left\{ (t,x): (\eta ^{(2)} _t (x) - \eta ^{(2)} _{t-} (x) ) >0  \right\}
\]
holds a.s. In other words, every moment of birth for $(\eta ^{(1)} _t)$ is a moment
of birth for $(\eta ^{(2)} _t)$ as well, and the spatial location of the birth is also identical.

\end{prop}

\textbf{Proof}. We can show by induction that each moment 
of birth for $(\eta ^{(1)} _t)$ is a moment
of birth for $(\eta ^{(2)} _t)$ as well,
and that each moment $\tau$ of death for $(\eta ^{(2)} _t)$
is a moment of death for $(\eta ^{(1)} _t)$ provided
$\eta ^{(1)} _{\tau -} (x) = \eta ^{(2)} _{\tau -} (x) $,
where $x$ is the site where the death at $\tau$ takes place.
Moreover, in both cases
the birth or the death occurs at the same site. Here a moment
of birth is a random time at which the value of 
the process at one of sites is increased by $1$,
and a
moment of death is a random time at which the value of 
the process at one of sites is decreased by $1$.
 The statement formulated
above implies \eqref{culprit}.

Denote by $\{ \tau _m \} _{m \in \N}$ the moments of 
 jumps of $(\eta ^{(1)}_t)$ and $(\eta ^{(2)}_t)$,
$0<\tau _1 <\tau _2 < \tau _3 < ...$. In other words, a time
$t\in \{ \tau _m \} _{m \in \N}$ if and only if 
at least one of the processes 
$(\eta ^{(1)}_t)$ and $(\eta ^{(2)}_t)$ jumps at the time $t$.

Here we deal only with
the base case, the induction step is done in the same way.
There is nothing to show if $\tau _1$ is a moment of birth
for $(\eta ^{(2)}_t)$ or a moment of death for $(\eta ^{(1)}_t)$.
Assume that $\tau _1$ is a moment of birth for $(\eta ^{(1)}_t)$
and let $x$ be the place of birth:
\[
 \eta ^{(1)} _{\tau _1} (x) - \eta ^{(1)} _{\tau _1 - } (x)  = 1.
\]
Note that
\[
 \eta ^{(1)} _{\tau _1 -} 
 = \eta ^{(1)} _{0} \leq \eta ^{(2)} _{0} =  \eta ^{(2)} _{\tau _1 -}.
\]
The process $(\eta ^{(1)}_t)$ satisfies \eqref{flippant}, hence
\[
 1 =\int\limits _{(0,\tau _1] \times [0, \infty ) }
I _{ [0, b _k(x,\eta _{s-} )] } (u) N_1 (x,ds,du) - 
\int\limits _{(0,\tau _1) \times [0, \infty ) }
I _{ [0, b _k(x,\eta _{s-} )] } (u) N_1 (x,ds,du)
\]
\[
  = \int\limits _{\{\tau _1\} \times [0, \infty ) }
I _{ [0, b _k(x,\eta _{s-} )] } (u) N_1 (x,ds,du)
\]
and
\[
 N_1(\{x\} \times \{\tau _1 \}\times [0, b_1 (x,\eta ^{(1)} _{0} )])=1 \ \ \ \ \text{ a.s. }
\]
Since $b_2  (x,\eta ^{(2)} _{0} ) \geq  b_1 (x,\eta ^{(1)} _{0} )$, 
$\tau _1 $ is a moment of birth at $x$ for $(\eta ^{(2)}_t)$.
The case when $\tau _1$ is a moment of death for $(\eta ^{(2)}_t)$
at a site $x$
and  $\eta ^{(2)} _{\tau _1 -} (x) = \eta ^{(2)} _{\tau _1 -} (x)$
is analyzed similarly.
\qed

 For $x \in \Z ^\d$, $\xi, \eta \in \X$ we define 
 \begin{equation} \label{tilde d}
   \tilde d(x, \xi, \eta) = \left\{
     \begin{array}{l l}
        d(x, \xi) \quad \text{if} \quad \xi (x) > \eta (x), \\
        d(x, \eta) \quad \text{if} \quad \xi (x) < \eta (x), \\
        d(x, \xi) \wedge d(x, \eta) \quad \text{if} \quad \xi (x) = \eta (x),
     \end{array} \right.
 \end{equation}
and
 \begin{equation} \label{tilde b}
   \tilde b(x, \xi, \eta) = \left\{
     \begin{array}{l l}
        b(x, \xi) \quad \text{if} \quad \xi (x) > \eta (x), \\
        b(x, \eta) \quad \text{if} \quad \xi (x) < \eta (x), \\
        b(x, \xi) \vee b(x, \eta) \quad \text{if} \quad \xi (x) = \eta (x),
     \end{array} \right.
 \end{equation}

 Note that $\tilde d(x, \xi, \eta) = \tilde d(x, \eta, \xi)$,
 $\tilde b(x, \xi, \eta) = \tilde b(x, \eta, \xi)$.
We will see below in \eqref{alacrity} and \eqref{alacrity2} how these
functions come into play.

\begin{lem} \label{schlemmen}
 For every $x \in \Z ^\d$, $\xi, \eta \in \X$,
 \begin{equation*}
  \tilde b(x, \xi, \eta) -  b(x, \eta) \leq \sum\limits _{y \in \Z ^\d} a(x-y) |\xi (x) -\eta (x)|,
 \end{equation*}
and
 \begin{equation*}
  \tilde d(x, \xi, \eta) -  d(x, \eta) \geq - \sum\limits _{y \in \Z ^\d} a(x-y) |\xi (x) -\eta (x)|.
 \end{equation*}
 
\end{lem}
\textbf{Proof}.
We have
\[
\tilde b(x, \xi, \eta) -  b(x, \eta) = I_{\{ \xi (x) > \eta (x) \}} (b(x, \xi) -  b(x, \eta))
+ I_{\{ \xi (x) = \eta (x) \}} \big[ (b(x, \xi) -  b(x, \eta)) \vee 0 \big]
\]
\[
 \leq \sum\limits _{x \in \Z ^\d} a(x-y) |\xi (x) -\eta (x)|.
\]
Similarly, 
\[
 \tilde d(x, \xi, \eta) -  d(x, \eta) = I_{\{ \xi (x) > \eta (x) \}} (d(x, \xi) -  d(x, \eta))
 +I_{\{ \xi (x) = \eta (x) \}} \big[ (d(x, \xi) -  d(x, \eta)) \wedge 0 \big]
\]
\[
 \geq - \sum\limits _{x \in \Z ^\d} a(x-y) |\xi (x) -\eta (x)|.
\]

 The next lemma will play the key role in the proof of 
 pathwise uniqueness for $\eqref{se lat}$.

 \begin{lem} \label{bloviate}
 Let $(( \xi _t )_{t\in [0,T]} , N_1 ,  N _2)$, $(\Omega , \mathscr{F} , P) $,
 $(\{ \mathscr {F} _t  \} _ {t\in [0,T]}) $ and 
 $(( \zeta _t )_{t\in [0,T]} , N_1 ,  N _2)$, $(\Omega , \bar {\mathscr{F}} , P) $,
 $(\{ \bar{\mathscr {F}} _t  \} _ {t\in [0,T]}) $
 be weak solutions to \eqref{se lat}. 
 Then 
 \begin{align} \label{bening}
 \begin{split}
 E & \sum\limits _{x \in \Z ^\d} w (x) |\xi _t (x) - \zeta _t (x)|  \\
  \leq 4 C_{w,a} \int\limits _{(0,t]  } 
   ds E
 \sum\limits_{x \in \Z ^\d}  w (x) & |\xi _{s-}(x) - \zeta _{s-}(x)| + 
 E  \sum\limits _{x \in \Z ^\d} w (x) |\xi _0 (x)  - \zeta _0 (x)|.
  \end{split}
 \end{align}
 
 \end{lem}

 \textbf{Proof}.
 Let
 $(\xi \vee _t \zeta ) $ be the c\`adl\`ag process defined by
 \[
  (\xi \vee _t \zeta) (x) = \xi _t (x) \vee \zeta_t (x),  \quad \quad t \in [0,T], \ x \in \Z ^\d.
 \]
 This process is adapted to the filtration $\{ \hat {\mathscr {F}} _t \}$,
 where $\hat {\mathscr {F}} _t : = \mathscr{F} _t \vee \bar{ \mathscr{F}} _t$.
 Note that since $N_1$ and $N_2$ are compatible with  $ \{\hat{\mathscr {F}} _t \} $
 	by Definition \ref{path uniq def}.

Define also
 \begin{equation} \label{alacrity}
  \tilde  d_t (x) := \int\limits _{(0,t] \times [0, \infty ) }
I _{ [0, \tilde d(x, \xi_{r-}, \zeta _{r-} )] } (v) N_2 (x,dr,dv)
 \end{equation}
 and
 \begin{equation}\label{alacrity2}
  \tilde  b_t (x) := \int\limits _{(0,t] \times [0, \infty ) }
I _{ [0, \tilde b(x, \xi_{s-}, \zeta _{s-} )] } (u) N_1 (x,ds,du).
 \end{equation}
 Then
 $\tilde b_t (x)$ and $\tilde d_t (x)$ are the numbers of births and deaths, respectively,
 for the process $(\xi \vee _t \zeta )$ at site $x$ that occurred before time $t$, 
 that is, 
  \begin{equation*}
 \tilde  d_t (x) = \# \left\{ r: r \leq t, \xi \vee _{r} \zeta  (x) - \xi \vee _{r-} \zeta  (x) = - 1 \right\},
  \end{equation*}
and similarly for $ \tilde b _t(x)$.
 Indeed, let $\tau$ be a moment of birth for $(\xi \vee _{t} \zeta )$,
 that is, $\xi \vee _{\tau} \zeta  (x) - \xi \vee _{\tau-} \zeta  (x) =  1$.
Without loss of generality assume that $\xi  _{\tau-} (x) \geq \zeta  _{\tau-} (x) $.
If $\xi  _{\tau-} (x) > \zeta  _{\tau-} (x) $, then $\tau$
is a moment of birth for $(\xi _t)$, hence 
$N_1(\{ x\} \times \{\tau \} \times [0,b(x, \xi _{\tau -})] ) = 1$ a.s.
and $\tilde b_{\tau} (x) -\tilde b_{\tau -} (x) =1 $.
If $\xi  _{\tau-} (x) = \zeta  _{\tau-} (x) $,
then $\tau$
is a moment of birth for at least one of the 
processes $(\xi _t)$ and  $(\zeta _t)$, hence 
\[
N_1(\{ x\} \times \{\tau \} \times [0,b(x, \xi _{\tau -}) \vee  b(x, \zeta _{\tau -}) ] ) = 1 \ \ \ \text{ a.s. }
\]
and again $\tilde b_{\tau} (x) - \tilde b_{\tau -} (x) =1 $.
On the other hand, let $\tilde b_{\tau} (x) - \tilde b_{\tau -} (x) = 1$.
Again, with no loss of generality we assume that 
$\xi  _{\tau-} (x) \geq \zeta  _{\tau-} (x) $.
If $\xi  _{\tau-} (x) > \zeta  _{\tau-} (x) $, then 
\[
{\tilde b (x, \xi  _{\tau-} (x), \zeta  _{\tau-} (x) )
= b (x, \xi  _{\tau-} )},
\]
hence
${N_1(\{ x\} \times \{\tau \} \times [0,b(x, \xi _{\tau -})] ) = 1}$ a.s.
and $\tau $ is a moment of birth for $(\xi \vee _t \zeta )$.
The remaining case $\xi  _{\tau-} (x) = \zeta  _{\tau-} (x) $
is similar.
The proof of \eqref{alacrity} follows the same pattern.

 Fix $t \in [0,T]$ and $x \in \Z ^\d$.
Note that 
 \[
 E\int\limits _{(0,t] \times [0, \infty ) } \big\{
I _{ [0,\tilde b(x,\xi _{s-} ,\zeta _{s-} )] }  (u) 
- I _{ [0, b(x,\zeta _{s-} )] } (u) \big\} N_1 (x,ds,du)
\]
 \[
= \int\limits _{(0,t]  } E
 I\{ \xi_{s-} (x) \geq \zeta_{s-} (x) \} \big\{
\tilde b(x,\xi _{s-} ,\zeta _{s-} ) 
- b(x,\zeta _{s-} ) \big\} ds
\]
 \[
 \leq  \int\limits _{(0,t]  }  ds
 E
 \sum\limits_{y \in \Z ^\d} a(x-y)|\xi _{s-}(y) - \zeta _{s-}(y)|,
 \]
 and
 \[
E \int\limits _{(0,t] \times [0, \infty ) } \big\{ 
I _{ [0, \tilde d(x,\xi _{s-} ,\zeta _{s-} )] } (v)-
I _{ [0,d(x,\zeta _{r-} )] } (v) \big\}N_2 (x,dr,dv) 
 \]
 \[
 = \int\limits _{(0,t]  } E
I\{ \xi_{s-} (x) \geq \zeta_{s-} (x) \}\big\{ 
 \tilde d(x,\xi _{s-} ,\zeta _{s-} ) -
d(x,\zeta _{r-} ) \big\} dr 
 \]
 \[
 \geq - \int\limits _{(0,t]  } 
   ds E
 \sum\limits_{y \in \Z ^\d} a(x-y)|\xi _{s-}(y) - \zeta _{s-}(y)|.
 \]
 So, we can write
 
 \[
  0 \leq E(\xi \vee _{t} \zeta  (x) - \zeta _t (x)) 
 \]
 \[
 = E\int\limits _{(0,t] \times [0, \infty ) } \big\{
I _{ [0,\tilde b(x,\xi _{s-} ,\zeta _{s-} )] }  (u) 
- I _{ [0, b(x,\zeta _{s-} )] } (u) \big\} N_1 (x,ds,du)
\]
\[
-E \int\limits _{(0,t] \times [0, \infty ) } \big\{ 
I _{ [0, \tilde d(x,\xi _{s-} ,\zeta _{s-} )] } (v)-
I _{ [0,d(x,\zeta _{s-} )] } (v) \big\}N_2 (x,ds,dv) 
+ E(\xi \vee _0 \zeta  (x) - \zeta _0 (x))
 \]
  \[
 \leq 2 \int\limits _{(0,t]  } 
   ds E
 \sum\limits_{y \in \Z ^\d} a(x-y)|\xi _{s-}(y) - \zeta _{s-}(y)|
 + E(\xi \vee _0 \zeta  (x) - \zeta _0 (x)).
 \]
 
 Multiplying the last inequality by $w (x)$
 and summing over $x$, we get

  \begin{align} \label{clandestine}
E \sum\limits _{x \in \Z ^\d} w (x) (\xi \vee _t \zeta  (x) - \zeta _t (x)) 
 \leq 2 & \int\limits _{(0,t]  } 
   ds E
 \sum\limits_{y \in \Z ^\d} |\xi _{s-}(y) - \zeta _{s-}(y)|
 \sum\limits _{x \in \Z ^\d} w (x) a(x-y)  \notag 
 \\
 + E  \sum\limits _{x \in \Z ^\d} w (x) |\xi \vee _0 \zeta (x)  - \zeta _0 (x)| 
 & \leq 2 \int\limits _{(0,t]  } 
   ds E
 \sum\limits_{y \in \Z ^\d} C_{w,a} w (y) |\xi _{s-}(y) - \zeta _{s-}(y)|
 \\
 +E & \sum\limits _{x \in \Z ^\d} w (x) |\xi \vee _0 \zeta  (x)  - \zeta _0 (x)|. \notag
 \end{align}

 Keeping in mind that $ (p \vee q -p) + (p \vee q - q) =|p-q|$, 
 we obtain \eqref{bening}
by swapping $\xi $ and $\zeta$ in \eqref{clandestine} and then adding
 the obtained inequality to \eqref{clandestine}. 
\qed

 \begin{prop} \label{un}
   Pathwise uniqueness holds for equation \eqref{se lat}.
 \end{prop}

 \textbf{Proof}.
 Let
 $(( \xi _t )_{t\in [0,T]}$
 and $(( \zeta _t )_{t\in [0,T]}$
 be two solutions to \eqref{se lat} 
 as in Lemma \ref{bloviate}.
 We know by
 item (ii) of Definition \ref{weak solution lat}
 that 
 $$
 f(t):={ E  \sum\limits _{x \in \Z ^\d} w (x) 
 |\xi _t (x) - \zeta _t (x)|} < \infty.
 $$
 Note that $f$ is a continuous function by the dominated convergence theorem, since for a fixed $s >0$
 every solution $(\eta _t)$ of \eqref{se lat}
 satisfies 
$\eta _{s-} = \eta _s = \eta _{s+} $ a.s.
 Furthermore, $f(0) = 0$, therefore 
 Grownwall's inequality and Lemma \ref{bloviate} yield $f(t)=0$.
 Since $\zeta _t (x) , \xi _t (x)$ are c\`adl\`ag processes, it follows that $\zeta _t (x) = \xi _t (x)$ a.s.
 for \textit{all} $t\in (0,T]$. 
 \qed
 
 Define $\bar b (x,\eta) := \sup\limits _{\alpha \leq \eta} b(x,\alpha) $.
Note that $\bar b$ is non-decreasing in the sense that
\begin{equation*}
 \bar b (x,\eta ^1) \leq \bar b (x,\eta ^2) \quad \text{ whenever } \eta ^1 \leq \eta ^2,
\end{equation*}
and that 
$\bar b$ satisfies inequalities of the form \eqref{bigot}.
Indeed, if $\xi, \eta \in \X, x \in \Z ^\d, \xi (x) \geq \eta (x)$, 
then
\[
 \bar b(x, \xi) - \bar b(x, \eta ) 
 =  \sup\limits _{\alpha : \ \alpha \leq \xi} 
 \left[ b(x,\alpha ) - \sup \limits _{ \beta : \ \beta \leq \eta } b(x,\beta ) \right]
 \]
 \[
 \leq 
 \sup\limits _{\alpha : \ \alpha \leq \xi } 
 \left[ b(x,\alpha ) -  b(x,\alpha \wedge  \eta ) \right]
 \leq \sup\limits _{\alpha: \  \alpha \leq \xi } 
\sum\limits _{y \in \Z ^\d} a(x-y) |\alpha (y) - \alpha (y) \wedge \eta (y) |
\]
\[
 \leq \sum\limits _{y \in \Z ^\d} a(x-y) |   \xi (y) - \eta (y) |.
\]

Also, for every $x\in \Z ^\d$
the map
\begin{align*}
 \X \ni \xi \mapsto \bar b(x, \xi)   \in \R _+ \notag
\end{align*}
is continuous by Lemma \ref{brisk},
since $\bar b$ satisfies all conditions
imposed on $b$.

 Before proceeding to the general existence result,
 let us consider a pure birth equation 
 \begin{equation}\label{pure birth se lat}
\begin{split}
\xi _t (x) = \int\limits _{(0,t] \times [0, \infty ) }
I _{ [0,\bar b(x,\xi _{s-} )] } (u) N_1 (x,ds,du) + \eta _0 (x).
\end{split}
\end{equation}
This equation is of the form \eqref{se lat}.

\begin{lem}
 Equation \eqref{pure birth se lat} has a (unique) solution.
\end{lem}

\textbf{Proof}.
Let us start with the
  equation \eqref{pure birth se lat} with a `truncated' initial condition
  and birth rate, that is, with the initial condition 
 \[
 \eta ^{(n)} _0 (x) = I_{\{|x|_1 \leq n \}} \eta  _0 (x) 
 \]
 and the birth rate 
 \[
 \bar b^{(n)}(x, \eta) = I_{\{|x|_1 \leq n \}}  \bar b(x, \eta).
 \]
 Here $n$ is a natural number.
 The initial condition is 
 finite and the birth rate satisfies
 \eqref{smolder},
 hence there exists a unique solution
 by Proposition \ref{finite in cond}. We denote this unique solution of 
 \begin{equation}\label{pure birth se lat (n)}
\begin{split}
\xi _t (x) = \int\limits _{(0,t] \times [0, \infty ) }
I _{ [0,\bar b ^{(n)} (x,\xi _{s-} )] } (u) N_1 (x,ds,du) + I_{\{|x|_1 \leq n \}} \eta _0 (x)  ,
\end{split}
\end{equation}
 by
 $(\xi ^{(n)} _t)_{t \in [0,T]}$. By Proposition \ref{comparison finite} we have
 $\xi ^{(m)} _t \leq \xi ^{(n)} _t$, $m \leq n$,
 and
\[
 \left\{ t : \sum\limits _{x \in \Z ^\d} (\xi ^{(m)} _t(x) -\xi ^{(m)} _{t-}(x) ) =1 \right\}
 \subset \left\{ t : \sum\limits _{x \in \Z ^\d} (\xi ^{(n)} _t(x) -\xi ^{(n)} _{t-}(x) ) =1 \right\}
\]
almost surely.
Therefore, the limit $\bar \eta _t = \lim\limits _{n \to \infty} \xi ^{(n)} _t$
exists and is c\`adl\`ag (if finite). For each $n\in \N$
\[
 E \sum\limits _{x \in \Z ^\d} w (x) \xi ^{(n)} _t(x) = 
 E \sum\limits _{x \in \Z ^\d} w (x) I_{\{|x|_1 \leq n \}} \int\limits _{(0,t] \times [0, \infty ) }
I _{ [0,\bar b (x,\xi ^{(n)} _{s-} )] } (u) N_1 (x,ds,du) 
\]
\[
+
E \sum\limits _{x \in \Z ^\d} w (x) I_{\{|x|_1 \leq n \}} \eta _0 (x)
\]
\[
 \leq  E\sum\limits _{x \in \Z ^\d} w (x) I_{\{|x|_1 \leq n \}} \int\limits _{(0,t]  }
\bar b  (x,\xi ^{(n)} _{s-} )  ds +
E \sum\limits _{x \in \Z ^\d} w (x) \eta _0 (x).
\]
Recall that $\0 \in \X$, $\0 (x) \equiv 0$.
By \eqref{bigot}
\[
\bar b(x,\xi ^{(n)} _{s-} ) \leq \sum\limits _{y \in \Z ^\d} 
a(x-y) \xi ^{(n)}(y) +  b(x,\0 ),
\]
hence
\[
  E\sum\limits _{x \in \Z ^\d} w (x) I_{\{|x|_1 \leq n \}} \int\limits _{(0,t]  }
\bar b  (x,\xi ^{(n)} _{s-} )  ds 
 \leq  E\sum\limits _{x \in \Z ^\d} w (x) \int\limits _{(0,t]  } ds
 \big[ \sum\limits _{y \in \Z ^\d} 
a(x-y) \xi ^{(n)} _{s-} (y) + b(x,\0 ) \big]
 \]
 \[
  \leq t \sum\limits _{x \in \Z ^\d} w (x) b(x,\0 )+ E \int\limits _{(0,t]  } ds \sum\limits _{y \in \Z ^\d} \xi ^{(n)} _{s-} (y)
  \sum\limits _{x \in \Z ^\d} w (x) a(x-y)
 \]
  \[
 \leq t \sum\limits _{x \in \Z ^\d} w (x) b(x,\0 ) +  C_{w , a} E \int\limits _{(0,t]  } ds 
 \sum\limits _{y \in \Z ^\d} w(y) \xi ^{(n)} _{s-}(y).
  \]
Thus, 
 \begin{equation}
  \begin{split}
   E & \sum\limits _{x \in \Z ^\d} w (x) \xi ^{(n)} _t (x) \\
   \leq  C_{w , a} E \int\limits _{(0,t]  } ds 
 \sum\limits _{x \in \Z ^\d} w(x) \xi ^{(n)} _{s-}(x) &+ t \sum\limits _{x \in \Z ^\d} w (x) b(x,\0 )+
   E \sum\limits _{x \in \Z ^\d} w (x) \eta _0 (x).
  \end{split}
 \end{equation}

The expression on the left hand side is finite by 
Proposition \ref{finite in cond} and depends continuously on $t$
by the same argument as in the proof of Proposition \ref{un}, therefore 
Grownwall's inequality implies
\begin{equation}
 E \sum\limits _{x \in \Z ^\d} w (x) \xi ^{(n)} _t (x)
 \leq e^{C_{w , a} t} \big[ E \sum\limits _{x \in \Z ^\d} w (x) \eta _0 (x) + 
 t \sum\limits _{x \in \Z ^\d} w (x) b(x,\0 ) \big].
\end{equation}
Letting $ n \to \infty$, we get by the monotone convergence theorem  
\begin{equation}
 E \sum\limits _{x \in \Z ^\d} w (x) \bar \eta _t (x)
 \leq e^{C_{w , a} t} \big[ E \sum\limits _{x \in \Z ^\d} w (x) \eta _0 (x) + 
 t \sum\limits _{x \in \Z ^\d} w (x) b(x,\0 ) \big].
\end{equation}
Since $b^{(n)}(x,\xi ^{(n)} _{s-}) \uparrow b(x,\bar \eta _{s-})$ a.s.,
$( \bar \eta _t)$ is a solution to
\eqref{pure birth se lat}. Uniqueness follows from Proposition \ref{un}.
 \qed

 \begin{prop} \label{ex}
  Strong existence holds for equation \eqref{se lat}.
 \end{prop}
\textbf{Proof}.
As in the proof of the previous proposition, we first consider 
  equation \eqref{se lat} with the 'truncated' initial condition
 $\eta ^{(n)} _0 (x) = I_{\{|x|_1 \leq n \}} \eta ^{(n)} _0 (x) $
 and the birth rate
  \[
  b^{(n)}(x, \eta) = I_{\{|x|_1 \leq n \}} b(x, \eta).
 \]

 We denote the unique solution of 
 \begin{equation}\label{se lat (n)}
\begin{split}
\eta _t (x) = \int\limits _{(0,t] \times [0, \infty ) }
I _{ [0,b ^{(n)} (x,\eta _{s-} )] } (u) N_1 (x,ds,du) \\
- \int\limits _{(0,t] \times [0, \infty ) }
I _{ [0,d(x,\eta _{r-} )] } (v) N_2 (x,dr,dv) + I_{\{|x|_1 \leq n \}} \eta _0 (x)  ,
\end{split}
\end{equation}
 by
 $(\eta ^{(n)} _t)_{t \in [0,T]}$.

 
 Let $m, n \in \N$. The estimations below are more natural when $m \leq n$,
 but formally we cover both cases.
 For $x \in \Z ^\d$  we have
 by \eqref{bigot}
 \begin{equation}\label{toad}
 \begin{gathered}
   E \sup\limits _{t \in [0,T]} \Big[
 \int\limits _{(0,t] \times [0, \infty ) } 
I _{ [ b ^{(m)} (x,\eta ^{(m)} _{s-} ),
b ^{(m)} (x,\eta ^{(m)} _{s-} ) \vee b ^{(n)} (x,\eta ^{(n)} _{s-} )] }  (u) 
 N_1 (x,ds,du) \Big] 
\\
  = E  
 \int\limits _{(0,T] \times [0, \infty ) } 
I _{ [ b^{(m)} (x,\eta ^{(m)} _{s-} ),
b ^{(m)}(x,\eta ^{(m)} _{s-} ) \vee b ^{(n)} (x,\eta ^{(n)} _{s-} )] }  (u) 
  N_1 (x,ds,du) 
\\
  = E 
 \int\limits _{(0,T]  } I\{ |x|_1 \leq m \} \big\{
b  (x,\eta ^{(m)} _{s-} ) \vee b  (x,\eta ^{(n)} _{s-} ) -b ^{(m)}(x,\eta ^{(m)} _{s-} )  \big\}
 ds 
\\
 +
 E 
 \int\limits _{(0,T]  } I\{ m< |x|_1 \leq n \} 
 b (x,\eta ^{(n)} _{s-} ) 
 ds 
\\
  \leq E
 \int\limits _{(0,T]  }  
\sum\limits_{y \in \Z ^\d} a(x-y)E|\eta ^{(n)} _{s-}(y) - \eta ^{(m)} _{s-}(y)| 
 ds 
\\
 + 
 E 
 \int\limits _{(0,T]  }  I\{ m< |x|_1 \leq n \} 
 \big\{
 b(x,\0) + \sum\limits_{y \in \Z ^\d} a(x-y)E\eta ^{(n)} _{s-}(y)
 \big\} ds.
 \end{gathered}
 \end{equation}
 On the other hand, as in the proof of Lemma \eqref{bloviate},
 
  \begin{equation}\label{toad2}
 \begin{gathered}
 E \inf\limits _{t \in [0,T]} \Big[ \int\limits _{(0,t] \times [0, \infty ) }  \big\{ 
I _{ [0, \tilde d(x,\eta ^{(n)} _{r-} ,\eta ^{(m)} _{r-} )] } (v)-
I _{ [0,d(x,\eta ^{(m)} _{r-} )] } (v) \big\}N_2 (x,dr,dv) \Big]
\\
  \geq
  - E \int\limits _{(0,T]  } 
\sum\limits_{y \in \Z ^\d} a(x-y)|\eta ^{(n)} _{r-}(y) - \eta ^{(m)} _{r-}(y)|
 dr.
 \end{gathered}
\end{equation}
 Therefore, by \eqref{alacrity}, \eqref{alacrity2}, \eqref{toad}, and \eqref{toad2},
  \begin{equation}\label{toad3}
\begin{gathered}
  E \sup\limits _{t \in [0,T]}  \big( \eta ^{(n)} \vee _t \eta ^{(m)}  (x) - \eta ^{(m)} _t (x) \big)
\\
  = E \sup\limits _{t \in [0,T]} \Big[
 \int\limits _{(0,t] \times [0, \infty ) } \big\{
I _{ [0, b^{(n)} (x,\eta ^{(n)} _{s-})  \vee b^{(m)}(x, \eta ^{(m)} _{s-} )] }  (u) 
- I _{ [0, b(x,\eta ^{(m)} _{s-} )] } (u) \big\} N_1 (x,ds,du)
\\
- \int\limits _{(0,t] \times [0, \infty ) } \big\{ 
I _{ [0, \tilde d(x,\eta ^{(n)} _{r-} ,\eta ^{(m)} _{r-} )] } (v)-
I _{ [0,d(x,\eta ^{(m)} _{r-} )] } (v) \big\}N_2 (x,dr,dv)
+  \big( \eta ^{(n)} \vee _0 \eta ^{(m)}  (x) - \eta ^{(m)} _0 (x) \big) \Big]
\\
\leq 2 \int\limits _{(0,T] \times [0, \infty ) } 
   ds E
 \sum\limits_{y \in \Z ^\d} a(x-y)|\eta ^{(m)} _{s-}(y) - \eta ^{(n)} _{s-}(y)| 
 + E \big( \eta ^{(n)} \vee _0 \eta ^{(m)}  (x) - \eta ^{(m)} _0 (x) \big)
\\
 +  E 
 \int\limits _{(0,T]  } ds  I\{ m< |x|_1 \leq n \} 
 \big\{
 b(x,\0) + \sum\limits_{y \in \Z ^\d} a(x-y)E\eta ^{(n)} _{s-}(y) .
 \big\}
 \end{gathered}
\end{equation}

By Proposition \ref{comparison finite}, a.s. $\eta ^{(n)} _{s-} \leq \xi ^{(n)} _{s-}$,
$s \geq 0$. 
Multiplying \eqref{toad3} by $w (x)$ and taking the sum over $x$, we obtain

 \[
 E  \sum\limits _{x \in \Z ^\d} 
  w (x) \sup\limits _{t \in [0,T]} \big( \eta ^{(n)} \vee _t \eta ^{(m)}  (x) - \eta ^{(m)} _t (x) \big)
 \]
\[
\leq 2 \int\limits _{(0,T] \times [0, \infty ) } 
   ds E \sum\limits _{x \in \Z ^\d} w (x)
 \sum\limits_{y \in \Z ^\d} a(x-y)|\eta ^{(m)} _{s-}(y) - \eta ^{(n)} _{s-}(y)| 
 + E \sum\limits _{x \in \Z ^\d}  w (x) \big( \eta ^{(n)} \vee _0 \eta ^{(m)}  (x) - \eta ^{(m)} _0 (x) \big)
\]
\[
 +  E \sum\limits _{x \in \Z ^\d} 
  w (x) I\{ m< |x|_1 \leq n \}
 \int\limits _{(0,T]  } ds
 \big\{
 b(x,\0) + \sum\limits_{y \in \Z ^\d} a(x-y)E\xi ^{(n)} _{s-}(y) 
 \big\}
\]
\[
\leq 2 \int\limits _{(0,T] \times [0, \infty ) } 
   ds E
 \sum\limits_{y \in \Z ^\d} w (y) C _{w , a} |\eta ^{(m)} _{s-}(y) - \eta ^{(n)} _{s-}(y)| 
 + E \sum\limits _{x \in \Z ^\d}  w (x) \big( \eta ^{(n)} \vee _0 \eta ^{(m)}  (x)- \eta ^{(m)} _0 (x) \big)
\]
\[
 + T b(x,\0) \sum\limits _{x \in \Z ^\d} 
  w (x) I\{ m< |x|_1 \leq n \} 
  +
  E  \int\limits _{(0,T]  } ds \sum\limits _{y \in \Z ^\d} \xi _{s-} (y) 
  \sum\limits _{x \in \Z ^\d} 
  w (x) a(x-y) I\{ m< |x|_1 \leq n \}.
\]

Using the above inequality
and the inequality 
\begin{equation*}
 \sup\limits _{t} |p _t - q _ t| \leq \sup\limits _{t} (p _t \vee q _ t - q _ t)
 +
 \sup\limits _{t} (p _t \vee q _ t - p _ t),
\end{equation*}
where $p, q$ are some functions with common domain and the supremum is taken over 
their domain,
we get
 
 \begin{align*} 
 \begin{split}
 \Delta _{m,n} :=E    \sum\limits _{x \in \Z ^\d} & w (x)
 \sup\limits _{t \in [0,T]} |\eta ^{(n)} _t (x) - \eta ^{(m)} _t (x)|  
 \\
 \leq 
 E  \sum\limits _{x \in \Z ^\d} 
  w (x) \sup\limits _{t \in [0,T]} \big( \eta ^{(n)} \vee _t \eta ^{(m)}  (x) - \eta ^{(m)} _t (x) & \big) + 
  E  \sum\limits _{x \in \Z ^\d} 
  w (x) \sup\limits _{t \in [0,T]} \big( \eta ^{(n)} \vee _t \eta ^{(m)}  (x) - \eta ^{(n)} _t (x) \big)
  \\
  \leq 4 C_{w,a} \int\limits _{(0,T] \times [0, \infty ) } 
   ds E
 \sum\limits_{x \in \Z ^\d}  w (x) |\eta ^{(n)} _{s-}(x) - & \eta ^{(m)} _{s-}(x)| 
 + E  \sum\limits _{x \in \Z ^\d}  w (x) |\eta ^{(n)} _0 (x) - \eta ^{(m)} _0 (x)| 
 \\
 + T b(x,\0) \sum\limits _{x \in \Z ^\d} 
  w (x) I\{ m< |x|_1 \leq n \} 
  + 
  E  \int\limits _{(0,T]  } & ds  \sum\limits _{y \in \Z ^\d}  \xi _{s-} (y) 
  \sum\limits _{x \in \Z ^\d} 
  w (x) a(x-y) I\{ m< |x|_1 \leq n \}
  \end{split}
 \end{align*}
 and
 consequently
   \begin{align} \label{pariah}
 \begin{split}
 \Delta _{m,n} 
  \leq \exp \{4 C_{w,a} T\}
  E \Bigg[ \sum\limits _{x \in \Z ^\d}  w (x) |\eta ^{(n)} _0 (x) - \eta ^{(m)} _0 (x)| +
  T b(x,\0) \sum\limits _{x \in \Z ^\d} 
  w (x) I\{ m< |x|_1 \leq n \} 
  \\
  + 
    \int\limits _{(0,T]  } ds \sum\limits _{y \in \Z ^\d}  \xi _{s-} (y) 
  \sum\limits _{x \in \Z ^\d} 
  w (x) a(y-x) I\{ m< |x|_1 \leq n \}
  \Bigg]
  \end{split}
 \end{align}
 by Gronwall's inequality.
 As $m,n \to \infty$, 
 $E \sum\limits _{x \in \Z ^\d}  w (x) |\eta ^{(n)} _0  - \eta ^{(m)} _0 (x)| \to 0$
 and
 $$\sum\limits _{x \in \Z ^\d} 
  w (x) I{\{ m< |x|_1 \leq n \}} \to 0.$$
  
  To deal with the third summand on the right hand side of \eqref{pariah},
  we define
   \[
   r(y, m) :=  \frac{\sum\limits _{x \in \Z ^\d} 
  w (x) a(x-y) I\{  |x|_1 > m \}}{ \sum\limits _{x \in \Z ^\d} 
  w (x) a(x-y) } .
  \] 
  
  Clearly,
  for each $y \in \Z ^\d$, $r(y, m) \to 0$ as $m \to \infty$.
  Hence 
  \[
   E  \int\limits _{(0,T]  } \sum\limits _{y \in \Z ^\d}  \xi _{s-} (y) 
  \sum\limits _{x \in \Z ^\d} 
  w (x) a(x-y) I\{ m< |x|_1 \leq n \} ds \leq 
 C_{w,a} E  \int\limits _{(0,T]  } \sum\limits _{y \in \Z ^\d} w(y) \xi _{s-} (y) r(y,m) ds
 \to 0
  \]
by \eqref{foible} and the dominated convergence theorem.

 Consequently,
   \begin{align} \label{saggy}
\Delta _{m,n}
\to 0, \quad  m,n \to \infty.
 \end{align}
 Since $w (x) >0$ for all $x$, \eqref{saggy} implies that
 \[
  P \{ \sup\limits _{t \in [0,T] } |\eta ^{(n)} _t (x) - \eta ^{(m)} _t (x)| >0 \} = 
  P \{ \sup\limits _{t \in [0,T] } |\eta ^{(n)} _t (x) - \eta ^{(m)} _t (x)| \geq 1 \}
 \]
\[
 \leq \frac{ \Delta _{m,n} }{w (x)} \to 0, \quad  m,n  \to \infty.
\]
Cauchy convergence in probability 
implies existence of a subsequence
  along which almost
 sure convergence takes place;
 moreover, using the diagonal argument, 
 we can find a subsequence $\{ n_m\} \subset \N$
 such that 
 for each $x$ there exists $(\eta _t (x))_{t \in [0,T]}$ satisfying
 \begin{equation} \label{verschwenderisch}
  P \{ \sup\limits _{t \in [0,T] } |\eta ^{(n_k)} _t (x) - \eta  _t (x)|  \to 0, k \to \infty \} = 1. 
 \end{equation}
Furthermore, 
$\eta ^{(n)} _t \leq \bar \eta  _t$, $t \in [0,T]$, where  $(\bar \eta  _t)$
is the unique solution of \eqref{pure birth se lat}.
Thus, since $\eta _t \leq \bar \eta  _t$, $t \in [0,T]$,
by the dominated convergence theorem

\begin{equation}
 P \{ \sup\limits _{t \in [0,T] }
 \sum\limits _{x \in \Z ^\d} w (x) |\eta ^{(n_k)} _t (x) - \eta  _t (x)| \to 0, k \to \infty \} =1. 
\end{equation}

 Since $b,d$ are continuous, $\eta ^{(n_{k})} _t \to \eta _t$ a.s. in $\X$
 and  $b^{(n_k)}(x,\eta ^{(n_k)} _{s-}) = b(x,\eta ^{(n_k)} _{s-})$ whenever $ n_k \geq |x|_1 $,
 $(\eta  _t)_{t \in [0,T]} $ is a strong solution to \eqref{se lat} if we can
 show that
 $E \sum\limits _{x \in \Z ^\d} w (x) \sup\limits _{t \in [0,T]} \eta  _t (x)  < \infty$, 
 the  integrals on the right hand side of \eqref{se lat} are well defined and 
  \begin{equation} \label{finite integrals}
  E \int\limits _0 ^T
  [b(x,\eta _{s-}) + d(x,\eta _{s-})] ds < \infty.
 \end{equation}
  The inequality
  $E \sum\limits _{x \in \Z ^\d} w (x) \sup\limits _{t \in [0,T]} \eta  _t (x)  < \infty$
  follows from the inequalities
  \[
  \eta ^{(n)} _t \leq \bar \eta _t, \quad n \in \N,
  \]
  where $(\bar \eta _t)$ is a solution to \eqref{pure birth se lat}.
  The integrals on the right hand side of \eqref{se lat} are well defined
  as pointwise limits of the corresponding integrals for $(\eta ^{(n_{k})} _t)$.

To prove \eqref{finite integrals}, we denote the number of births and deaths at $x$
before $t$ by $b_t (x)$ and $d_t (x)$ respectively, i.e.
\[
 b_t (x) = \# \{ s: \eta _{s}(x) - \eta _{s-}(x) = 1 \}=
\int\limits _{(0,t] \times [0, \infty ) }
I _{ [0,b(x,\eta _{s-} )] } (u) N_1 (x,ds,du).
\]
and similarly for $d_t (x)$. Note that 
$\eta _{t} (x) = b_t (x) - d _t(x) + \eta _0 (x)$.
Let $(\tau _n )$ be the moments of jumps of 
  $c_t (x) := b_t (x) + d _t(x)$,
  $0 = \tau _0 \leq \tau _1 \leq \tau _2 \leq ... $.
  If $\tau _k <T$
  and $c_t (x) = c_{\tau _k} (x) $ for all $t \in [\tau _k ,T]$, 
 we set $\tau _{k+j} =  T$ for all $j \in \N$.
  Note that $\tau _n$ is a stopping time
  with respect to the filtration $\{ \mathscr{S}_t \}$.
 We have
  \[
   c_t(x) = \sum\limits _{n \in \N} I \{ \tau _n \leq t \}
  \]
  a.s. for all $t \in [0,T)$.
 Define for $n \in \N$
\begin{equation*}
\begin{split}
c ^{(n)} _t (x) := \int\limits _{(0,t] \times [0, \infty ) }
I _{ [0,b(x,\eta _{s-} ) \wedge n] } (u) N_1 (x,ds,du) \\
+ \int\limits _{(0,t] \times [0, \infty ) }
I _{ [0,d(x,\eta _{r-} ) \wedge n] } (v) N_2 (x,dr,dv) .
\end{split}
\end{equation*}
  Then 
  \[
   M^{(n)}_t (x) = c ^{(n)} _t (x) - \int\limits _0 ^t 
\big(b(x,\eta _{s-} ) \wedge n \big) ds - \int\limits _0 ^t 
\big(d(x,\eta _{s-} ) \wedge n \big) ds
  \]
is a martingale with respect to
$\{ \mathscr{S} _t \}$. By the optional stopping theorem
$EM^{(n)}_{\tau _1 } (x) = 0 $, 
hence
\[
 E \int\limits _0 ^{\tau _1 }
\big( b(x,\eta _{s-} ) \wedge n + d(x,\eta _{s-} ) \wedge n \big) ds  \leq 1.
\]
Similarly,
\[
 E \int\limits _{\tau _{m}} ^{\tau _{m+1} }
\big( b(x,\eta _{s-} ) \wedge n + d(x,\eta _{s-} ) \wedge n \big) ds \leq P\{\tau _{m} < T  \}.
\]
Consequently
\[
 E \int\limits _{0} ^{T }
\big( b(x,\eta _{s-} ) \wedge n + d(x,\eta _{s-} ) \wedge n \big) ds 
\leq
\sum\limits _{m = 0} ^ \infty E \int\limits _{\tau _{m}} ^{\tau _{m+1} }
\big( b(x,\eta _{s-} ) \wedge n + d(x,\eta _{s-} ) \wedge n \big) ds
\]
\[
 \leq \sum\limits _{m = 0} ^ \infty  P\{\tau _{m} < T  \} \overset{*}{=} 
 \sum\limits _{m = 0} ^ \infty  P\{c_T (x) \geq m  \} = E c_T (x) + 1,
\]
where the transition marked by the asterisk is possible
in particular since
\[
 c_T (x) =  c_{T-} (x) \quad \text{ a.s.}
\]

Letting $n \to \infty$, we get by the monotone convergence theorem
\begin{equation} \label{cleavage}
 E \int\limits _{0} ^{T }
\big( b(x,\eta _{s-} ) + d(x,\eta _{s-} )  \big) ds
\leq E c_T (x) + 1
\end{equation}

Since only existing particles may disappear, the number of deaths $d_t (x)$
satisfies for every $t \in [0,T]$
\[
 d_t (x) \leq b_t (x) + \eta _0 (x).
\]
Finally, since by Proposition \ref{comparison finite} every birth for $(\eta_t)$ is a birth at the same time and place for $(\bar \eta_t)$
as well (note that Proposition \ref{comparison finite} cannot be applied  to $(\eta_t)$ and $(\bar \eta_t)$ directly, but 
 can be to the processes  $(\eta ^{(n)} _t)_{t \geq 0}$
 and  $(\xi ^{(n)} _t)_{t \geq 0}$), we have a.s. $b_T (x) \leq \bar \eta _T (x)$, and hence
\begin{equation} \label{on the skids}
 E c_T (x) \leq 2 E b_T (x) + E \eta _0 (x)
\leq 
2 E \bar \eta _T (x) + E \eta _0 (x) < \infty.
\end{equation}
\qed
  
  \begin{rmk} \label{livid}
   In fact,  \eqref{cleavage} and \eqref{on the skids}
   yield even stronger inequality 
   \[
    \sum\limits _{x \in \Z ^\d} w (x) E \int\limits _{0} ^{T }
\big( b(x,\eta _{s-} ) + d(x,\eta _{s-} )  \big) ds < \infty.
   \]

  \end{rmk}

 The following statement is a consequence of 
  Proposition \ref{un} and 
 \cite[Theorem 3.14]{KurtzYamWatan}
 
 \begin{prop} \label{uniqueness in law thm}
  Joint uniqueness in law holds for \eqref{se lat}.
 \end{prop}

 \begin{prop}\label{Markov property}
  The unique solution to
  \eqref{se lat} is
  a Markov process:
  for all $ \mathscr{D} \in D _{\mathcal{X}} [0,\infty) $
  and $q\geq 0$,
   \begin{equation} \label{fizzle out}
    P \big[ (\eta _{q+\cdot}) \in \mathscr{D} \mid \mathscr{S} _q  \big]= 
    P \big[ (\eta _{q+\cdot}) \in \mathscr{D} \mid \eta _q  \big].
   \end{equation}
  
 \end{prop}
 
\textbf{Proof}. For $t \geq q $ we have 

 \begin{equation}
\begin{split}
\eta _t (x) = \int\limits _{(q,t] \times [0, \infty ) }
I _{ [0,b(x,\eta _{s-} )] } (u) N_1 (x,ds,du) \\
- \int\limits _{(q,t] \times [0, \infty ) }
I _{ [0,d(x,\eta _{r-} )] } (v) N_2 (x,dr,dv) + \eta _q (x),
\end{split}
\end{equation}
 therefore $(\eta _{q+\cdot})$ is 
 $\sigma \{ \eta_q , N_k(\{x\} \times [q,q+r] \times C), 
x\in \Z ^\d, C\in \mathscr{B} (\R_+), r \geq 0 , k=1,2 \}$-measurable
by Propositions \ref{un} and \ref{ex}; the fact that we start from the time 
$q$ instead of $0$ does not cause problems. 
Since Poisson processes have independent increments, 
\[
\sigma \{ N_k(\{x\} \times [q,q+r] \times C), 
x\in \Z ^\d, C\in \mathscr{B} (\R_+), r \geq 0 , k=1,2 \}
\]
is independent of $\mathscr{S} _q$ and
 \eqref{fizzle out} follows. \qed

\section{Proof of Propositions \ref{poignant} and \ref{filibuster}}

  \textbf{Proof of Proposition} \ref{poignant}.
  For $R>0$ we define
$\mathbf{B}_R := \{ x \in \Z ^\d \mid |x| _1 \leq R \}$.
  By Ito's formula
  
 \begin{equation} \label{engender}
  \begin{split}
F (\eta _t) = F (\eta _0) + 
\int\limits _{\mathbf{B}_{R_F} \times (0,t] \times [0, \infty ) } \big\{ F(\eta ^{+x} _{s-}) - F(\eta _{s-}) \big\}
I _{ [0,b(x,\eta _{s-} )] } (u) N_1 (dx,ds,du) \\
  + \int\limits _{ \mathbf{B}_{R_F} \times (0,t] \times [0, \infty ) } \big\{ F(\eta ^{-x} _{r-}) - F(\eta _{r-}) \big\}
I _{ [0,d(x,\eta _{r-} )] } (v) N_2 (dx,dr,dv).
\end{split}
 \end{equation}
We can write 
\[
\int\limits _{\mathbf{B}_{R_F} \times (0,t] \times [0, \infty ) } \big\{ F(\eta ^{+x} _{s-}) - F(\eta _{s-}) \big\}
I _{ [0,b(x,\eta _{s-} )] } (u) N_1 (dx,ds,du)
\]
\[
= 
\int\limits _{(0,t] } \sum\limits _{x \in \mathbf{B}_{R_F}} 
\big\{ F(\eta ^{+x} _{s-}) - F(\eta _{s-}) \big\}
b(x,\eta _{s-} ) ds
\]
\[+ \int\limits _{\mathbf{B}_{R_F} \times (0,t] \times [0, \infty ) }
\big\{ F(\eta ^{+x} _{s-}) - F(\eta _{s-}) \big\}
I _{ [0,b(x,\eta _{s-} )] } (u) \tilde N_1 (dx,ds,du)
\]
where
$\tilde N_1 = N_1 - \# (dx) ds du$.
Since $F(\eta ^{+x} ) - F(\eta )$ is bounded uniformly
 in $x$ and $\eta $, the last integral with respect to $\tilde N_1$
 is a martingale by item (iv) of Definition \ref{weak solution lat},
 see e.g. \cite[Section 3 of Chapter 2]{IkedaWat}. 
Similarly,
\[
\int\limits _{ \mathbf{B}_{R_F} \times (0,t] \times [0, \infty ) } 
\big\{ F(\eta ^{-x} _{r-}) - F(\eta _{r-}) \big\}
I _{ [0,d(x,\eta _{r-} )] } (v) N_2 (dx,dr,dv)
 \]
 can be represented as a sum of 
 \[
\int\limits _{  (0,t]  }
\sum\limits _{x \in \mathbf{B}_{R_F}}
\big\{ F(\eta ^{-x} _{r-}) - F(\eta _{r-}) \big\}
d(x,\eta _{r-} ) dr
 \]
 and a martingale. The assertion of the proposition now follows
 from \eqref{engender} and \eqref{the generator}.

To prove Proposition \ref{filibuster}
 we will need the following form of 
 the martingale representation theorem, which is a corollary of
 \cite[Theorem 7.4, Chapter 2]{IkedaWat}.
 
 \textbf{Theorem}. \emph{Let $(\alpha _t)$ 
 be an increasing c\`adl\`ag $\mathcal{X}$-valued process on a filtered probability space 
  $(\Omega, \mathscr{F},{ \{ \mathscr{F}_t}, t\geq 0 \}, P)$
 such that the point Process  defined by 
 \[
 Q_p([0,t]\times \{x\}) = \alpha _t (x), \ \ \ t \geq 0,  \ \ \ x \in \Z ^\d
 \]
  has the (predictable) 
  compensator $p (t, E) = \int\limits_{x \in E, s \in [0,t]}  \phi (x,s) \#(dx)  ds$
 such that
 \[
 E p (t,\{x\}) < \infty
 \] 
 for each $t\geq 0$ and $x \in \Z ^\d$.
 Furthermore, assume that a.s. there are no simultaneous jumps:
 $${\sup\limits _{t \geq 0} \sum\limits _{x \in \Z ^\d}[\alpha _t(x) - \alpha _{t-}(x)] \leq 1}.$$
 Then on an extended filtered probability space
 $(\tilde \Omega, \tilde{ \mathscr {F}}, 
 \{ \tilde{\mathscr {F}}_t, t\geq 0 \}, \tilde P)$
 there exists an adapted to $\{ \tilde{\mathscr {F} }_t, t\geq 0 \} $ 
 Poisson point process $\tilde N$
 on $\Z ^\d \times \R _+ \times \R _+$ 
  with 
 intensity measure $\# \times ds \times du $
 such that
 \begin{equation}
  \alpha _t (x) = \int\limits _{ (0,t] \times \R _+ }
  I_{[0,\phi(x,s)]} (u) \tilde N (x, ds,du), \ \ \ x \in \Z ^\d.
 \end{equation}}
 To see that this theorem follows from 
 \cite[Theorem 7.4, Chapter 2]{IkedaWat},
 we should take there
 $\mathbf{X} = \Z ^\d$, $\mathbf{Z} = \Z ^\d \times \R _+$,
 $m = \# \times du$, $q(t,E) = \phi(t,E)  $,
 $\theta (t, (x,u)) ={ x I \{ u \leq \phi (x,t) \}} + 
 {\Delta I \{ u > \phi (x,t) \}}$.

 \textbf{Proof of Proposition} \ref{filibuster}. 
 Define
 a $(\Z _+)^{\Z ^\d \times \{-1,1 \}}$-valued
 process
 $ \{ \alpha _{t}, t\geq 0 \}$ 
 by
\begin{align*}
 \alpha _t ((x,1)) = \eta ^{(b)} _{t} (x), \quad
 \alpha _t ((x,-1)) = \eta ^{(d)} _{t} (x).
\end{align*}
Conditions of the previous theorem are satisfied, 
so we get
\[
\alpha _t ((x,1)) = \int\limits _{ (0,t] \times \R _+ }
  I_{[0,b(x, \eta _{s-})]} (u)  N ((x,1), ds,du),
\]
and 
\[
\alpha _t ((x,-1)) = \int\limits _{ (0,t] \times \R _+ }
  I_{[0,d(x, \eta _{s-})]} (u)  N ((x,-1), ds,du),
\]
a.s. for all $t \in [0,T]$,
where $N$ is a Poisson point process on $(\Z ^\d \times \{-1,1 \}) \times \R _+ \times \R _+$.
Define $N_1$ and $N_2$ by 
\[
 N_1 (x \times [0,t] \times U)= N ((x,1)\times [0,t] \times U),
\]
\[
 N_2 (x \times [0,t] \times U)= N ((x,-1)\times [0,t] \times U).
\]
Then a.s. for all $t \in [0,T]$
\[
 \eta _t (x) = \eta _0 (x) + \eta ^{(b)} _{t} (x) - \eta ^{(d)} _{t} (x)
\]
\[
 = \int\limits _{(0,t] \times [0, \infty ) }
I _{ [0,b(x,\eta _{s-} )] } (u) N_1 (x,ds,du)
- \int\limits _{(0,t] \times [0, \infty ) }
I _{ [0,d(x,\eta _{r-} )] } (v) N_2 (x,dr,dv).
\]
 \qed

  \section{Proof of Theorem \ref{schmooze}
  and Proposition \ref{beseech}}\label{inv measure}

 Let us recall that $(\eta (\alpha, t))_{t \in [0,T]}$
 is the unique solution to \eqref{se lat} 
 with initial condition 
$\eta _0 \equiv \alpha$, $\alpha \in \X$.

 \begin{lem} \label{martingaleV}
  The process 
   \[
 M_t := V (\eta (\0, t)) - \int\limits _0 ^t L V (\eta (\0, s-))ds,
\]
 is well defined and an $\{ \mathscr{F}_t \}$-martingale.
 \end{lem}

 \textbf{Proof}. 
 Denote
 \[
  \mathcal{D} (L) := \{ \eta \in \X : 
  \sum\limits _{x \in \Z ^\d} v (x) [b(x , \eta) + d (x,\eta)] < \infty \}.
 \]

 For $\eta \in  \mathcal{D} (L)$ the expression 
 $L V (\eta)$ in \eqref{the generator}
 is well defined. 
 Since $v$ satisfies 
 the same assumptions as 
 $w$ does, Theorem \ref{core thm},
  Proposition  \ref{poignant} and 
 all the other results proven in Sections 3 and 4 are still valid
 if we replace in their formulations 
 $w$ by $v$ and $\X$ by $\X _v$.
Remark \ref{livid}
 implies that a.s.
 \[
  E \int\limits _0 ^T L V (\eta (\0, s-))ds < \infty,
 \]
 in particular, $ \eta (\0, s) \in \mathcal{D} (L)$ a.s.
 since $(\eta (\0, t))$ is c\`adl\`ag.
 By Proposition \ref{poignant}, 
\[
 \sum\limits _{x \in \Z ^\d} v (x) \eta (\0, t) (x) I \{ |x|_1 \leq n \} - 
 \int\limits _0 ^t \sum\limits _{x \in \Z ^\d}  v (x) 
 \big\{ b(x,\eta (\0, s-)) - d(x ,\eta (\0, s-)) \big\} I \{ |x|_1 \leq n \} ds
\]
is an $\{ \mathscr{F}_t \}$-martingale. By the dominated convergence theorem,
\[
 \sum\limits _{x \in \Z ^\d} v (x)\eta (\0, t) (x) I \{ |x|_1 \leq n \}
 \overset{L^1}{ \to} \sum\limits _{x \in \Z ^\d} v (x) \eta (\0, t) (x).
\]
Furthermore,
\[
 \int\limits _0 ^t \sum\limits _{x \in \Z ^\d}  v (x) 
 \big\{ b(x,\eta (\0, s-)) - d(x ,\eta (\0, s-)) \big\} I \{ |x|_1 \leq n \}  ds
\]
 \[
  \overset{L^1}{ \to}
   \int\limits _0 ^t \sum\limits _{x \in \Z ^\d}  v (x) 
 \big\{ b(x,\eta (\0, s-)) - d(x ,\eta (\0, s-)) \big\} ds
\]
 since the difference goes to zero in $L^1$ by Remark \ref{livid}.
 Therefore, 
 \[
 M_t = V (\eta (\0, t)) - \int\limits _0 ^t L V (\eta (\0, s-))ds,
\]
 is an $\{ \mathscr{F}_t \}$-martingale. 
 \qed

\textbf{Proof of Theorem} \ref{schmooze}. 
 For $\alpha \in \X$  let us define
 \[
   P^t (\alpha, B) := P \{ \eta (\alpha, t) \in B \}, 
   \ \ \ B \in \mathscr{B} (\X), t \geq 0
 \]
and let
 $$
 \mu _ n (B) := \frac 1n \int\limits _0 ^n P ^s \{ \0, B \} ds, \quad B \in \mathscr{B}(\X).
 $$
 Denote also $K_r :=\{\eta \in \X _v : V (\eta) \leq r \}$, $r>0$.
 Imitating the proof of Lemma 9.7 of Chapter 4 \cite{EthierKurtz}, we obtain
 by Lemma \ref{martingaleV}
 
 \[
 0 \leq E V (\eta _n) = E V ( \eta _0) + E \int ^n _0 LV (\eta _{s-}) ds
 \]
 \[
  =E V ( \eta _0) + E \int ^n _0 LV (\eta _{s-}) I \{\eta _{s-} \in K_r \} ds + 
  E \int ^n _0 LV (\eta _{s-}) I \{\eta _{s-} \notin K_r \} ds 
 \]
\[
 \leq E V ( \eta _0) + E c_1 \int ^n _0  I \{\eta _{s-} \in K_r \} ds 
 + (c_1 - c_2 r) E \int ^n _0 I \{\eta _{s-} \notin K_r \} ds 
 \]
 \[
 = E V ( \eta _0) + n c_1 \mu _n ( K_r) + n (c_1 - c_2 r)[1- \mu _n ( K_r)],
\]
 hence
\begin{equation*}
 \mu _n ( K_r) \geq 1 - \frac{c_1}{c_2 r} - \frac{E V ( \eta _0)}{nc_2 r}.
\end{equation*}
We see that $\mu _n ( K_r) \to 1$ as $r \to \infty$ uniformly in $n \in \N$.
It follows from \eqref{tipsy}
that for every $r>0$ the set $K_r$ is precompact in $\X$,
therefore the family $\{\mu _n ,  n \in \N \}$ is tight.
 By Prohorov's theorem there exists 
a measure $\mu$ on $\X$
 and a sequence $\{n_k\}$ such that 
 $\mu _{n _k} \Rightarrow \mu$.
 Without loss of generality we assume that $\mu _n \Rightarrow \mu$.
 Let us show that 
 $\mu$ is an invariant measure.
 Take $f \in C _b (\X)$, then
 \[
  \int P ^t f(\eta) \mu (d\eta)
  =  \lim _n \int P ^t f(\eta) \mu _n (d\eta) =
  \lim _n \frac{1}{n}\int _0 ^{n} ds 
\int P ^t f(\eta) P ^s (\0, d \eta) =
 \]
 \[
 \lim _n \frac{1}{n}\int _0 ^{n} ds 
 P ^{t+s} f(\0) = 
\lim _n \frac{1}{n}\int _t ^{n+t} ds 
 P ^{s} f(\0)
\]
 \[
 = \lim _n \Big[ \frac{1}{n} \int _0 ^{n} + \frac{1}{n} \int _n ^{n+t} 
 - \frac{1}{n} \int _0 ^{t} \Big] = 
 \lim _n \int f(\eta) \mu _n (d \eta) = 
 \int  f(\eta) \mu  (d \eta).
\]
 \qed

 \textbf{Proof of Proposition} \ref{beseech}. 
 Let us take $v(x) = \frac{1}{1+|x|^{d+1} _1}$,
 and let $\mathbf{o} _\d$ be the origin in $\Z ^\d$.
 In the computations below
 we set ${C_1 = \sum\limits _{x \in \Z ^\d} \frac{b_0}{1+|x|^{d+1} _1}} $.
 Since $a_+$ satisfies a finite range property and $\sup\limits_{\substack{x, y \in \Z^\d:\\
 |x-y| <R}} \frac{v(x)}{v(y)} < \infty$, 
 there exists $C_2 >0 $
 such that 
 \[
 \sum\limits _{x \in \Z ^\d} v(x) a_+ (x-y) \leq C_2 v(y),  \ \ \ y \in \Z ^\d.
 \]
 Let also $m = m_1 \wedge m_2 >0$.
 We have
 for all $\eta \in \X _v$ 
 \[
  L V (\eta) \leq \sum\limits _{x \in \Z ^\d} v(x)[b_0 + \sum\limits _{y \in \Z ^\d} a_+ (x-y) \eta (y)]
  - m \sum\limits _{x \in \Z ^\d} v(x) \eta ^2 (x)
  \]
\[
\leq C_1 + \sum\limits _{y \in \Z ^\d}  \eta (y) \sum\limits _{x \in \Z ^\d} v(x) a_+ (x-y)
-
m \sum\limits _{x \in \Z ^\d} v(x) \eta ^2 (x) 
\leq C_1 + C_2 \sum\limits _{y \in \Z ^\d} v(y) \eta (y) 
-
m \sum\limits _{x \in \Z ^\d} v(x) \eta ^2 (x) 
\]
\[
\leq  C_3 - C _4 \sum\limits _{x \in \Z ^\d} v(x) \eta  (x) 
\]
for some constants $C_3, C_4 >0$.
Thus the desired statement follows 
from Theorem \ref{schmooze}.

 \section{Extinction and critical value for a model
 with independent branching birth rate and local death rate
 }\label{extinct surv}
 
In this section we consider the birth and death 
rates given by 
\begin{equation} \label{resurgence}
 b_{\lambda}(x,\eta ) = \lambda \sum\limits _{y: |y-x|\leq 1} \eta (y),
 \quad d(x,\eta ) = g(\eta (x)),
\end{equation}
where
$\lambda > 0$ and $g: \Z _+ \to \R_+$
  is a non-decreasing function 
  such that $g(0)=0$, $g(1)=1$  and $g(n)\geq n$.
  For  $g(n) = n ^2$ the evolution
  of 
 the process can be described as follows.
 Each particle is deleted from the system at a rate which
 is equal to the number of particles
 at the same site. Each particle gives
 birth to a new particle at a constant rate.
 The descendant appears at a site
 chosen uniformly among those
 neighboring to the predecessor sites and the site 
 of the predecessor.
 We denote the unique solution 
 of \eqref{se lat} by $(\eta ^{\lambda} _t) _{t \in [0,\infty)}$,
 or simply $(\eta _t) $.

 Let us consider equation
 \begin{equation}\label{se contact}
\begin{split}
\xi _t (x) = \int\limits _{(0,t] \times [0, \infty ) }
I _{ [0,b _{cont}(x,\xi _{s-} )] } (u) N_1 (x,ds,du) \\
- \int\limits _{(0,t] \times [0, \infty ) }
I _{ [0,d(x,\xi _{r-} )] } (v) N_2 (x,dr,dv) + \xi _0 (x)  ,
\end{split}
\end{equation}
 where $\xi _0 (x) = \eta _0 (x) \wedge 1$
 and $b _{cont}$ is given in \eqref{smug}.
 Equation \eqref{se contact} is of the form
 \eqref{se lat}. The unique solution
 $(\xi ^\lambda _t)_{t\in[0,T]}$
 of \eqref{se contact} is in fact
 the contact process, see Remark \ref{spin systems}.

 \begin{prop} \label{1224} Let $\lambda < \bar \lambda$. Then
 
 (i) $\xi ^\lambda _t \leq \eta ^\lambda _t$ a.s. for all $t \geq 0$,
 
 (ii) $\eta ^\lambda _t \leq \eta ^{\bar \lambda} _t   $ a.s. for all $t \geq 0$.
 \end{prop}

 \textbf{Proof}. We saw in
 the proof of Proposition \ref{ex}
 that every solution
 is an a.s. limit of solutions
 with finite initial conditions.
 Therefore,
 this statement is a consequence
 of Proposition \ref{comparison finite}. \qed
 
    The idea to couple the process with rates similar to
  \eqref{resurgence} with the contact process appeared in 
  Section 6.2
  \cite{FournierMeleard},
  however the rigorous proof has not been carried out there.

  We recall that $\mathbf{o}_\d$ stands for
  the origin in $\Z ^\d$.
  Let $\eta _0 (x) = I_{\{ x=\mathbf{o}_\d \}}$,
  and define

  \begin{equation*}
    p_s (\lambda) = P\{ \eta ^\lambda _t \ne \0 \text{ for all } t\geq 0 \}.
  \end{equation*}

From Proposition \ref{1224} it follows that
 $p_s$ is a non-decreasing function of $\lambda$. A
 standard comparison with a subcritical branching process 
 shows
 that $p_s (\lambda) = 0$ for sufficiently small $\lambda$, 
 for example for $\lambda < \frac{1}{2d+1}$. On the other hand,
 comparison with the contact process demonstrates that 
 $p_s (\lambda) >0 $ for $\lambda > \lambda _c ^{cont}$,
 where
 $\lambda _c ^{cont}$ is a critical value of the contact process.
 Therefore, there exists a critical value:
 
 \begin{equation*}
   \lambda _c = \inf \{ \lambda >0: p_s(\lambda) >0 \}.
 \end{equation*}

 We summarize the above discussion in the following proposition.
 
 \begin{prop}
  Consider the unique solution to \eqref{se lat} with the birth and 
  death rates \eqref{resurgence} 
  and the initial condition $\eta _0 (x) = I_{\{ x=\mathbf{o}_\d \}}$. 
  Then there exists $\lambda _c > 0$ such that
  
  (i) the process goes extinct if $\lambda < \lambda _c$:
  \[
  P\{ \eta _t = \0 \text{ for some } t\geq 0 \} =1,
  \]
  (ii) the process survives with positive probability if $\lambda > \lambda _c$:
  \[
  P\{ \eta _t \ne \0 \text{ for all } t\geq 0 \} >0.
  \]

 \end{prop}

 \section*{Acknowledgement}
 
 The authors acknowledge the financial support of the DFG through the SFB 701 
 ``Spektrale Strukturen
und Topologische Methoden in der Mathematik'' (Bielefeld University). 
V.B.
is also thankful the the support of the DFG through
  the IRTG (IGK) 1132 ``Stochastics and Real World Models''.
 V.B. and Y.K. would like to thank 
  Errico
 Presutti
 for the discussions during their visit to L'Aquila.

\bibliographystyle{alpha}
\bibliography{Sinus}

\newcommand{\etalchar}[1]{$^{#1}$}
\begin{thebibliography}{BRASS07}

\bibitem[AN72]{Branch1}
K.~B. Athreya and P.~E. Ney.
\newblock {\em Branching processes}.
\newblock Springer-Verlag, New York-Heidelberg, 1972.
\newblock Die Grundlehren der mathematischen Wissenschaften, Band 196.

\bibitem[And82]{And82}
E.~D. Andjel.
\newblock Invariant measures for the zero range processes.
\newblock {\em Ann. Probab.}, 10(3):525--547, 1982.

\bibitem[Bez14]{Bezborodovthesis}
V.~Bezborodov.
\newblock {\em Spatial birth-and-death {M}arkov processes}.
\newblock PhD thesis, Bielefeld University, 2014.

\bibitem[Bez15]{Bez15}
V.~Bezborodov.
\newblock Spatial birth-and-death markov dynamics of finite particle systems.
\newblock 2015.
\newblock http://arxiv.org/abs/1507.05804.

\bibitem[BFKR10]{Assumeexistence}
M.~Bal{\'a}zs, G.~Farkas, P.~Kov{\'a}cs, and A.~R{\'a}kos.
\newblock Random walk of second class particles in product shock measures.
\newblock {\em J. Stat. Phys.}, 139(2):252--279, 2010.

\bibitem[BP97]{BolkPacala}
B.~M. Bolker and S.~W. Pacala.
\newblock Using moment equations to understand stochastically driven spatial
  pattern formation in ecological systems.
\newblock {\em Theoretical Population Biology}, 52:179–197, 1997.

\bibitem[BP99]{BolkPacala2}
B.~M. Bolker and S.~W. Pacala.
\newblock Spatial moment equations for plant competitions: Understanding
  spatial strategies and the advantages of short dispersal.
\newblock {\em The American Naturalist}, 153(6):575–602, 1999.

\bibitem[BRASS07]{BRASS07}
M.~Bal{\'a}zs, F.~Rassoul-Agha, T.~Sepp{\"a}l{\"a}inen, and S.~Sethuraman.
\newblock Existence of the zero range process and a deposition model with
  superlinear growth rates.
\newblock {\em Ann. Probab.}, 35(4):1201--1249, 2007.

\bibitem[CT85]{misanthrope85}
C.~Cocozza-Thivent.
\newblock Processus des misanthropes.
\newblock {\em Z. Wahrsch. Verw. Gebiete}, 70(4):509--523, 1985.

\bibitem[DL05]{DieckmannLaw}
U.~Dieckmann and R.~Law.
\newblock Relaxation projections and the method of moments.
\newblock page 412–455. Cambridge University Press, 2005.

\bibitem[dLF06]{Yuleproc}
A.~de~La~Fortelle.
\newblock Yule process sample path asymptotics.
\newblock {\em Electron. Comm. Probab.}, 11:193--199 (electronic), 2006.

\bibitem[EH05]{EH05}
M.~R. Evans and T.~Hanney.
\newblock Nonequilibrium statistical mechanics of the zero-range process and
  related models.
\newblock {\em J. Phys. A}, 38(19):R195--R240, 2005.

\bibitem[EK86]{EthierKurtz}
S.~N. Ethier and T.~G. Kurtz.
\newblock {\em Markov processes}.
\newblock Wiley Series in Probability and Mathematical Statistics: Probability
  and Mathematical Statistics. John Wiley \& Sons, Inc., New York, 1986.
\newblock Characterization and convergence.

\bibitem[EK14]{EK14}
A.~M. Etheridge and T.~G. Kurtz.
\newblock Genealogical constructions of population models.
\newblock 2014.
\newblock preprint; arXiv:1402.6724 [math.PR].

\bibitem[Eth04]{Eth04}
A.~M. Etheridge.
\newblock Survival and extinction in a locally regulated population.
\newblock {\em Ann. Appl. Probab.}, 14(1):188--214, 2004.

\bibitem[FKK09]{FKK09}
D.~Finkelshtein, Y.~G. Kondratiev, and O.~Kutoviy.
\newblock Individual based model with competition in spatial ecology.
\newblock {\em SIAM J. Math. Anal.}, 41(1):297--317, 2009.

\bibitem[FKK12]{Semigroupapproach}
D.~Finkilstein, Y.~G. Kondratiev, and O.~Kutoviy.
\newblock Semigroup approach to birth-and-death stochastic dynamics in
  continuum.
\newblock {\em J. Funct. Anal.}, 262(3):1274–1308, 2012.

\bibitem[FKK13]{FKKBPDL}
D.~Finkelshtein, Y.~G. Kondratiev, and O.~Kutoviy.
\newblock An operator approach to {V}lasov scaling for some models of spatial
  ecology.
\newblock {\em Methods Funct. Anal. Topology}, 19(2):108--126, 2013.

\bibitem[FKKZ14]{Aggreg}
D.~Finkelshtein, Y.~G. Kondratiev, O.~Kutoviy, and E.~Zhizhina.
\newblock On an aggregation in birth-and-death stochastic dynamics.
\newblock {\em Nonlinearity}, 27(6):1105–1133, 2014.

\bibitem[FM04]{FournierMeleard}
N.~Fournier and S.~M{\'e}l{\'e}ard.
\newblock A microscopic probabilistic description of a locally regulated
  population and macroscopic approximations.
\newblock {\em Ann. Appl. Probab.}, 14(4):1880--1919, 2004.

\bibitem[GK06]{GarciaKurtz}
N.~L. Garcia and T.~G. Kurtz.
\newblock Spatial birth and death processes as solutions of stochastic
  equations.
\newblock {\em ALEA Lat. Am. J. Probab. Math. Stat.}, 1:281--303, 2006.

\bibitem[Hol70]{Hol70}
R.~Holley.
\newblock A class of interactions in an infinite particle system.
\newblock {\em Advances in Math.}, 5:291--309 (1970), 1970.

\bibitem[HS76]{HolStr76}
R.~A. Holley and D.~W. Stroock.
\newblock A martingale approach to infinite systems of interacting processes.
\newblock {\em Ann. Probability}, 4(2):195--228, 1976.

\bibitem[IW81]{IkedaWat}
N.~Ikeda and S.~Watanabe.
\newblock {\em Stochastic Differential Equations and Diffusion Processes}.
\newblock Nord-Holland publiching company, 1981.

\bibitem[Kal02]{Kallenbergfound}
O.~Kallenberg.
\newblock {\em Foundations of modern probability}.
\newblock Probability and its Applications. Springer-Verlag, second edition,
  2002.

\bibitem[KK02]{KondKuna}
Y.~G. Kondratiev and T.~Kuna.
\newblock Harmonic analysis on configuration space. {I}. {G}eneral theory.
\newblock {\em Infin. Dimens. Anal. Quantum Probab. Relat. Top.},
  5(2):201--233, 2002.

\bibitem[KKM10]{KKM10}
Y.~G. Kondratiev, O.~Kutoviy, and R.~Minlos.
\newblock Ergodicity of non-equilibrium {G}lauber dynamics in continuum.
\newblock {\em J. Funct. Anal.}, 258(9):3097--3116, 2010.

\bibitem[KKZ06]{KKZh06}
Y.~G. Kondratiev, O.~Kutoviy, and E.~Zhizhina.
\newblock Nonequilibrium {G}lauber-type dynamics in continuum.
\newblock {\em J. Math. Phys.}, 47(11):113501, 17, 2006.

\bibitem[KP96]{KP96}
T.~G. Kurtz and P.~E. Protter.
\newblock Weak convergence of stochastic integrals and differential equations.
  {II}. {I}nfinite-dimensional case.
\newblock In {\em Probabilistic models for nonlinear partial differential
  equations ({M}ontecatini {T}erme, 1995)}, volume 1627 of {\em Lecture Notes
  in Math.}, pages 197--285. Springer, Berlin, 1996.

\bibitem[KS05]{KesS05}
H.~Kesten and V.~Sidoravicius.
\newblock The spread of a rumor or infection in a moving population.
\newblock {\em Ann. Probab.}, 33(6):2402--2462, 2005.

\bibitem[KS08]{KesS08}
H.~Kesten and V.~Sidoravicius.
\newblock A shape theorem for the spread of an infection.
\newblock {\em Ann. of Math. (2)}, 167(3):701--766, 2008.

\bibitem[Kur07]{KurtzYamWatan}
T.~G. Kurtz.
\newblock The {Yamada-Watanabe-Engelbert} theorem for general stochastic
  equations and inequalities.
\newblock {\em Electron. J. Probab}, 12:951–965, 2007.

\bibitem[Lig85]{Lig85}
T.~M. Liggett.
\newblock {\em Interacting particle systems}.
\newblock Grundlehren der Mathematischen Wissenschaften. Springer, 1985.

\bibitem[LS81]{LS81}
T.~M. Liggett and F.~Spitzer.
\newblock Ergodic theorems for coupled random walks and other systems with
  locally interacting components.
\newblock {\em Z. Wahrsch. Verw. Gebiete}, 56(4):443--468, 1981.

\bibitem[MDL04]{MDL04}
D.~J. Murrell, U.~Dieckmann, and R.~Law.
\newblock On moment closures for population dynamics in continuous space.
\newblock {\em J. Theor. Biol.}, 229:421--432, 2004.

\bibitem[OFK{\etalchar{+}}14]{Ecology}
O.~Ovaskainen, D.~Finkelshtein, O.~Kutovyi, S.~Cornell, B.~M. Bolker, and Y.~G.
  Kondratiev.
\newblock A mathematical framework for the analysis of spatiotemporal point
  processes.
\newblock {\em Theoretical Ecology}, 7:101--113, 2014.

\bibitem[Pen08]{Pen08}
M.~D. Penrose.
\newblock Existence and spatial limit theorems for lattice and continuum
  particle systems.
\newblock {\em Probab. Surv.}, 5:1--36, 2008.

\bibitem[Spi70]{Spit70}
F.~Spitzer.
\newblock Interaction of {M}arkov processes.
\newblock {\em Advances in Math.}, 5:246--290 (1970), 1970.

\end{thebibliography}

\end{document}